\documentclass{article}
\usepackage{amsmath,amssymb,amsthm}
\theoremstyle{plain}
\newtheorem{thm}{Theorem}
\newtheorem{lm}{Lemma}
\newtheorem{prop}{Proposition}

\usepackage{caption,graphicx}
\usepackage[all]{xy}

\title{Analytic knots, satellites and the 4-ball genus}
\author{Burglind J\"oricke \\ { \ }}
\date{}

\begin{document}

\maketitle

\vglue 1cm

\begin{abstract}
\noindent Call a smooth knot (or smooth link) in the unit sphere in
$\mathbb{C}^2$ analytic (respectively, smoothly analytic) if it
bounds a complex curve (respectively, a smooth complex curve) in the
complex ball. Let $K$ be a smoothly analytic knot. For a small
tubular neighbourhood of $K$ we give a sharp lower bound for the
4-ball genus of analytic links $L$ contained in it.

\end{abstract}

\vfill

\noindent {\em Mathematics Subject Classification (2000)}

\noindent {\em Primary:} 57M25, 57M12, 32B15

\noindent {\em Secondary:} 14H30

\smallskip


\newpage

\section{Introduction}

Let $K$ be a smooth knot in the 3-sphere $S^3$. A smooth knot or
smooth link $L$ contained in a tubular neighbourhood $N(K) \subset
S^3$ of $K$ is a satellite of $K$ if it is not isotopic to $K$ in
$N(K)$ and (any connected component of $L$) is
not contained in a 3-ball inside $N(K).\,$ 
Satellites of non-trivial knots have been considered since long.
They play a role in the problem which knot complements admit
hyperbolic structure and received recent interest from the point of
view of invariants of knots and links. In the following we always
consider knots and links to be smooth and oriented. Speaking about
tubular neighbourhoods we will always consider them smoothly
bounded. Let $K$ be a non-trivial knot and $L$ a link contained in a
smoothly bounded tubular neighbourhood $N(K)$ of $K$.
Define an entire number $n$ in the following way. Consider a
projection ${\rm pr} : N(K) \to K$. The image ${\rm pr} (L) \subset
K$ is homologuous to $n \cdot [K]$ (in $H_1(K)$) for an entire
number $n$. Write $w_{N(K)} (L) = n$ and call $n$ the winding number
of $L$ in $N(K)$. (Having in mind a specified tubular neighbourhood
we will also speak about the winding number of $L$ around $K$.)

\smallskip

The pattern $\mathcal{L}$ of $L$ gives more precise information on
the satellite $L$. It is defined as follows. Denote by $U$ a
standard realization of the unknot. For instance, identify $S^3$
with the unit sphere $\partial {\mathbb B}^2$ in ${\mathbb C}^2$.
Let $H = \{ z_2 = 0 \}$ be the first coordinate line in ${\mathbb
C}^2$ and let $U$ be the unknot $U =
\partial {\mathbb B}^2 \cap H$ oriented as boundary of a complex
disc in $\mathbb{B}^2$. Consider a tubular neighbourhood $N (U)
\subset S^3$ of $U$. Trivialise both, $N(K)$ and $N(U)$, by Seifert
framing, i.e. by a transversal vector field on the knot which points
in the direction of a smooth oriented surface which is contained in
the sphere and bounded by the knot. Such a surface is called a
Seifert surface. Note that the trivialization does not depend on the
choice of the Seifert surface. Consider a diffeomorphism $\varphi_K
: N(K) \to N(U)$ which maps Seifert framing to Seifert framing. The
pattern $\mathcal{L}$ of the satellite $L$ is the isotopy class of
$\varphi_K (L)$ in $N(U)$.

\smallskip

A classical paper of Schubert \cite{S} relates the genus of a knot
to the genus of its satellites. The (smooth) genus $g(K)$ of a knot
(or link) $K$ is the minimal genus among smooth oriented surfaces in
$S^3$ bounded by $K$. (If $L$ is a link and the surface is not
connected we mean the sum of the genera of the connected
components.) The genus depends only on the isotopy class of the knot
or link. Define the genus of an isotopy class of links as the genus
of its representatives. Schubert's theorem is the following.

\bigskip

{\em For a satellite knot $L$ in a tubular neighbourhood $N(K)$ of a
knot $K \subset S^3$ with $n = w_{N(K)} (L)$ the following
inequality for the genera holds

\begin{equation}\label{eq1}
g(L) \geq \vert n \vert \, g \, (K) \, .
\end{equation}

Moreover,}

\begin{equation}\label{eq2}
g(L) \ge \vert n \vert \, g \, (K) + g(\mathcal{L}) \, .
\end{equation}

\bigskip

Identify again $S^3$ with the unit sphere $\partial {\mathbb B}^2$
in ${\mathbb C}^2$. We consider links (or knots) which are obtained
as the transverse intersection of $\partial {\mathbb B}^2$ with a
relatively closed complex curve $\tilde X$ in a neighbourhood of the
closed unit ball $\overline{{\mathbb B}^2}$. Following Rudolph
\cite{R} we call such links analytic, and smoothly analytic if the
complex curve $\tilde X \cap {\mathbb B}^2$ bounded by the link is
smooth (i.e., non-singular). We always consider an analytic link in
the sphere $\partial {\mathbb B}^2$ oriented as boundary of a
complex curve in the unit ball.

\smallskip

Since $H^2 ({\mathbb B}^2 , {\mathbb Z}) = 0$ the complex curve is the zero
locus $\{ z \in {\mathbb B}^2 : f(z) = 0 \}$ of an analytic function in a neighbourhood of the closed ball.
The function $f$ can be uniformly approximated on $\overline{{\mathbb B}^2}$ by a polynomial, which gives
an isotopic link in $\partial {\mathbb B}^2$ that bounds a
piece of an algebraic hypersurface. If the curve $\{ z \in {\mathbb B}^2 : f(z) = 0 \}$ is singular
its genus is defined to be the genus of its smooth perturbation $\{ z \in {\mathbb B}^2 : f(z) = \varepsilon \}$ for
generic small enough numbers $\varepsilon$.
\smallskip

We are interested in the (smooth) $4$-ball genus $g_4 (L)$ of a knot
(or link) $L$, called also slice genus. This is the minimal genus
among smooth oriented surfaces embedded into ${\mathbb B}^2$ and
bounded by $L$. Always $g_4 (L) \leq g(L)$ but $g_4 (L)$ may be
strictly smaller than $g(L)$. The $4$-ball genus gives a lower bound
for the unknotting number of a knot, the smallest number of crossing
changes needed to unknot the knot. The class of analytic knots is
interesting from the point of view of knot invariants: for them half
the Rasmussen invariant and also the $\tau$-invariant are equal to
the 4-ball genus of the knot (\cite{P}, \cite{He}, \cite{Sh}).

\smallskip

By a consequence of a deep theorem of Kronheimer
 and Mrowka
(Corollary 1.3 of \cite{KM}, the local Thom Conjecture) the 4-ball
genus of an analytic knot is realized on the complex curve bounded
by it. The proof of Kronheimer and Mrowka also shows that for a link
which bounds a connected complex curve in $\mathbb{B}^2$ its
$4$-ball genus is realized by the genus of this curve.

\smallskip

The following theorem holds.
\begin{thm}
Let $K$ be a smoothly analytic knot in $\partial {\mathbb B}^2$.
There exists a tubular neighbourhood $N(K) \subset \partial {\mathbb
B}^2$ of $K$ such that for any analytic link $L \subset N(K)$ the
number $n = w_{N(K)} (L)$ is
non-negative and the following statements hold.\\
\begin{enumerate}

\item[1.]  If $L$ is itself a knot then
\begin{equation}\label{eq4}
g_4 (L) \geq ng_4 (K) - \left[ \frac{n-1}{2}\, \right]\,.
\end{equation}
($\,[x]$ denotes the largest integer not exceeding the real number $x$).
\item[2.] Let $L$ be a link which bounds a connected complex curve $Y$. If $n$ is positive then the following lower bound for the $4$-ball
genus holds
\begin{equation}\label{eq3}
g_4 (L) \geq n g_4 (K) - (n-1) \, .
\end{equation}

\end{enumerate}

\noindent For $n=1$ the statements are true also if $K$ bounds a
singular curve.

\smallskip

\noindent The estimates are sharp in the following sense.

\begin{enumerate}

\item[3.] For each smoothly analytic knot $K$ with $g_4(K) \ge 1$, each natural number $n \ge 1$
and any tubular neighbourhood $N(K) \subset \partial \mathbb{B}^2 $
of $K$ there exists a link $L  \subset N(K)$ with $w_{N(K)}(L)=n$
which bounds a connected complex curve such that equality in
\eqref{eq3} is attained.

\item[4.] Further, for each smoothly analytic knot $K$ with $g_4(K) \ge 1$ and each natural number $n \ge 1$
there is a smoothly analytic knot $K_1$ which is smoothly isotopic
to $K$ and has the following property.  For any tubular
neighbourhood $N(K_1) \subset
\partial \mathbb{B}^2 $ of $K_1$ there exists an analytic knot $L \subset N(K_1)$ with winding number $w_{N(K)}(L)=n$ such that
equality in \eqref{eq4} is attained for the knot $K_1$ and the knot
$L$. In general, the original knot $K$ does not have this property.
\end{enumerate}

\end{thm}

The condition of analyticity of $K$ and $L$ cannot be removed.
Indeed, for any knot $K$ (in particular, for an analytic knot $K$)
and any tubular neighbourhood $N(K)$ the connected sum $L$ with its
mirror can be realized as a satellite of $K$ with $n = w_{N(K)} (L)
= 1$ and $g_4 (L) = 0$.

We do not know how big the tubular neighbourhood $N(K)$ in the
theorem can be chosen, in particular, we do not know whether for
$n=1$ the statement is true for analytic links in {\em any} tubular
neighbourhood of an analytic knot. For $n>1$ we do not know sharp
estimates of the $4$-ball genus of satellites of knots which bound
singular complex curves.

Some satellites are especially simple and useful. They are defined
in terms of closed braids. Recall the following definitions. Let
$C_n(\mathbb{C})= \{ (z_1, \ldots, z_n) \in \mathbb{C}^n:\, z_i \ne
z_j \; \rm{for} \; i \ne j\}\,$ be the configuration space of $n$
particles which move in the plane without collision. The symmetrized
configuration space $C_n(\mathbb{C}) \diagup \mathcal{S}_n$ is the
quotient of $C_n(\mathbb{C})$ by the action of the symmetric group
$\mathcal{S}_n$. Each point in $C_n(\mathbb{C}) \diagup
\mathcal{S}_n$ can be be considered as unordered tuple of $n$ points
and can be identified with the monic polynomial whose collection of
zeros equals this unordered tuple. The space of monic polynomials of
degree $n$ without multiple zeros is denoted by $\mathfrak{P}_n$,
the space of all monic polynomials of degree $n$ is denoted by
$\overline{\mathfrak{P}_n}$. The set of coefficients of polynomials
in $\mathfrak{P}_n$ is equal to $\mathbb{C}^n \setminus \{\sf{D}_n
=0 \}$, where $\sf{D}_n$ is the discriminant, i.e. $\sf{D}_n$  is a
polynomial on $\mathbb{C}^n$ which vanishes exactly if the monic
polynomial with these coefficients has multiple zeros.

Recall that a geometric braid with base point $E_n \in
\mathfrak{P}_n$ can be considered as a continuous map of the
interval $[0,1]$ into $\mathfrak{P}_n$ with initial and terminating
point equal to $E_n$ . A braid with base point $E_n$ is an isotopy
class of geometric braids with this base point. Such braids form a
group which is isomorphic to a group $\mathcal{B}_n$ with $n-1$
generators, denoted by $\sigma_1, \dots, \sigma_{n-1},\,$ and
finitely many relations. The group is called Artin's braid group.

A braid is quasi-positive if it is the product of conjugates of the
standard generators $\sigma_i$ of ${\mathcal B}_n$ (conjugates of
inverses of these generators are not allowed as factors).

We call an oriented closed curve ${\tilde L}$ in  $S^1 \times D^2$ a
closed geometric braid if the projection to $S^1$ is orientation
preserving on ${\tilde L}$. The circle $S^1$ is assumed to be
oriented, $D^2$ is a disc of real dimension $2$. The number of
preimages of a point under the projection is called the number of
strands. A closed braid is a free isotopy class of closed geometric
braids. (Free isotopy means isotopy without fixing a base point.) It
is well-known that free isotopy classes of closed geometric braids
on $n$ strands (for short, closed $n$-braids) are in one to one
correspondence to conjugacy classes in Artin's braid group
${\mathcal B}_n$ of braids on $n$ strands (see e.g. \cite{Bi}).

The notion of closed geometric braids is sometimes used in a more
special situation, namely, by a closed geometric braid one means an
oriented closed curve ${\tilde L}$ in $\partial {\mathbb B}^2
\backslash \{ z_1 = 0 \}$ for which $d \arg z_1 \mid {\tilde L} >
0$.

Note that alternatively a geometric braid with base point $E_n$ can
be considered as a collection of $n$ disjoint arcs in $[0,1] \times
D^2$ which join the collection $\{1\} \times E_n$ in the top $\{1\}
\times D^2$ with the "{identical}" {collection} $\{0\} \times E_n$
in the bottom $\{0\} \times D^2$ and is such that for each arc the
canonical projection to $[0,1]$ is a homeomorphism. Identifying top
and bottom we obtain a closed geometric braid in $S^1 \times D^2$,
called the closure of the geometric braid.

S. Orevkov pointed out that for $n > 1$ the statement of Theorem 1
does not extend to the situation of two closures of quasi-positive
geometric braids, one being a satellite of the other \cite{O1}. (For
convenience of the reader details are given below in example 3.)

Let a tubular neighbourhood $N(K)$ of the knot $K$ be the image of a
diffeomorphism from $S^1 \times D^2$ onto $N(K)$ so that $K$ is the
image of $S^1 \times \{0\}$. (We assume that the canonical framing
on $S^1 \times D^2$ is mapped to Seifert framing.) If a link $L$ in
$N(K)$ is the image of a closed geometric braid on $n$ strands in
$S^1 \times D^2$ then $L$ is called an $n$-braided link around $K$.

\smallskip


For an analytic link $L \subset \partial {\mathbb B}^2$ which bounds
a complex curve $Y$ and a $4$-ball $\Omega' \subset \partial
{\mathbb B}^2$ whose boundary is piecewise smooth and intersects $Y$
generically we call the intersection $Y \cap \partial \Omega' $ the
$\Omega'$-truncation of
$L$.

The following theorem describes all links $L$ which may appear in
the si\-tuation of Theorem~1. Let as above $K \subset \partial
{\mathbb B}^2$ be an analytic knot, $K= \tilde X \cap \partial
\mathbb{B}^2$ for a smooth relatively closed complex curve $\tilde
X$ in a neighbourhood of $\bar {\mathbb{B}}^2$, and let $L \subset
\partial {\mathbb B}^2$ be an analytic link, $L=\tilde Y \cap \partial {\mathbb B}^2$
for a relatively closed complex curve $\tilde Y$ in a neighbourhood
of $\bar {\mathbb{B}}^2$. Denote by $\mathbb D$ the unit disc in the
complex plane.

\begin{thm}
Let $K$ be a smoothly analytic knot in $\partial {\mathbb B}^2$.
There exists a tubular neighbourhood $N(K) \subset \partial
{\mathbb{B}}^2\;$ of $K$ and a pseudoconvex ball $\Omega' \subset
{\mathbb{B}}^2$ with piecewise smooth boundary which is obtained
from $\mathbb{B}^2$ by replacing a tubular neighbourhood of $K$ in
$\partial \mathbb{B}^2$ by a
Levi-flat hypersurface, 
such that 
for any analytic link $L\subset N(K)$
the following holds.
\begin{enumerate}
\item[1.] {\bf (Truncation.) }
The $\Omega'$-truncation $K' \subset
\partial \Omega'$ of $K$ is a smooth knot of the same $4$-ball genus $g_4(K') = g_4(K)$ as
$K$. If $n= w_{N(K)}(L)>0$ then the $\Omega'$-truncation $L' \subset
\partial \Omega'$ of $L$ is an $n$-braided link around
$K'$. If $n= w_{N(K)}(L)=0$, then $L'$ is the empty set.
The statements are true for $\Omega'$ replaced by a smoothly bounded
strictly pseudoconvex domain $\Omega_1,\; \Omega' \subset \Omega_1
\subset \mathbb{B}^2,\,$  (depending on $K$ and $L$),  with $C^2$ boundary which is
(away from corners of $\partial \Omega'$) $C^2$ close to
$\partial{\Omega}'$.



\item[2.] {\bf (Patterns of analytic closed braids in $N(K_1)$.)}
For a strictly pseudoconvex domain $\Omega_1$ as in statement 1
the pattern $\mathcal{L}_1$ of the $\Omega_1$-truncation $L_1$ of $L$ 
corresponds to the conjugacy class in the braid group ${\mathcal
B}_n$ of the product of
a quasi-positive braid $w \in {\mathcal B}_n$
with $g = g_4 (K)$ commutators in ${\mathcal B}_n$.

\end{enumerate}

\noindent The statements are sharp in the following sense.

\begin{enumerate}

\item[3.] {\bf (Realization of patterns as analytic links.)} Let $\mathcal{L}$ be a pattern
as described in statement 2. Then for each analytic knot there
exists an isotopic analytic knot $K \subset \partial \mathbb{B}^2$
such that the following holds. For any a priori given tubular
neighbourhood of $K$, the pattern $\mathcal{L}$ 
can be realized by an analytic link contained in this neighbourhood.
\end{enumerate}
\end{thm}

Notice that there is a piecewise smooth homeomorphism between
$\partial \Omega'$  and $\partial \mathbb{B}^2$, and $\partial
\Omega_1$ is diffeomorphic to $\partial \mathbb{B}^2$. The pattern
of links $L' \subset
\partial \Omega'$, and of links  $L_1 \subset
\partial \Omega_1$ repectively, contained in tubular neighbourhoods of knots
$K'$, and $K_1$ respectively, are defined using the piecewise smooth
homeomorphism, and the smooth homeomorphism respectively.

There is a continuous decreasing family $\Omega_t$, $t \in [0,1]$,
of strictly pseudoconvex balls $\Omega_t$ and a continuous family
$\psi_t$ of contactomorphisms $\psi_t : \partial \Omega_t \to
\partial {\mathbb B}^2$ such that $\Omega_0 = {\mathbb B}^2$ ,
$\psi_0 = {\rm id}$ and $\psi_t (K) = K_t \overset{\rm def}{=} X
\cap \partial \Omega_t$. So, one can ``identify'' the link $L_1$ of
statement 2 with a closed $n$-braid in a tubular neighbourhood of
$K$ in $\partial {\mathbb B}^2$ (for short, with an $n$-braided link
around $K$ in $\partial {\mathbb B}^2$ ).

\smallskip


The pattern $\mathcal{L}$ in statement 2 is not necessarily
quasi-positive but all quasipositive patterns can be realized in
statement 3. In the quasipositive case we have the following precise
statement on the $4$-ball genus of the satellite. 

\begin{lm}
Let $K \subset \partial {\mathbb B}^2$ be an analytic knot. Suppose
$L$ is a knot contained in a tubular neighbourhood $N(K)$ of $K$
whose pattern $\mathcal{L}$ is the closure of a quasi-positive
$n$-braid. Then 
$$
g_4 (L) = n \, g_4 (K) + g_4(\mathcal{L}) \, .
$$
\end{lm}



The results of this paper grew out of an unsuccessful attempt to
answer the question below. This question concerns the case of
strictly pseudoconvex domains instead of the ball $\mathbb{B}^2$,
and is related to the following fact (\cite{J}).\\
{\it Let $\Omega$ be a strictly pseudoconvex domain in a
two-dimensional Stein manifold. Then each element $e$ of the
fundamental group of the boundary $\pi_1 (\partial \Omega)$ whose
representatives are contractible in $\Omega$ can be represented by
the boundary of an immersed analytic disc in $\Omega$.}

For an immersed analytic disc with simple transverse
self-intersections the self-intersection number is the number of
double points. The number of self-intersections of a general
analytic disc is the number of self-intersections of its small
generic perturbations.
\bigskip

\noindent {\bf Question.} 1. What is the minimal self-intersection
number among analytic discs whose boundary represents $e$?

\noindent 2. In particular, let $\Omega = \{ (x,y,z) \in {\mathbb
C}^3 : x^2 + y^3 + z^5 = \varepsilon \} \cap {\mathbb B}^3$ be the
natural Stein filling of the Poincar\'e sphere. ($\varepsilon > 0$
is a small positive number, ${\mathbb B}^3$ is the unit ball in
${\mathbb C}^3$.) The loops $\{ x = \varepsilon^{\frac{1}{2}} \}
\cap \partial \Omega$, $\{ y = \varepsilon^{\frac{1}{3}} \} \cap
\partial \Omega$, $\{ z = \varepsilon^{\frac{1}{5}} \} \cap \partial
\Omega$ bound analytic discs in $\Omega$. Do these discs minimize
the self-intersection number among analytic discs whose boundaries
represent the respective element of $\pi_1 (\partial \Omega)$?

\bigskip

Proposition 1 below implies the following. Consider the very
restrictive class of analytic discs, whose boundaries are contained
in a small tubular neighbourhood of one of the loops in part 2 of
the question and are homotopic in $\partial \Omega$ to the
respective loop. Among them the mentioned discs have the minimizing
property of the self-intersection number. We do not know how to get
rid of the very restrictive condition.

Notice that the question concerns the complex structure of $\Omega$
rather than its Stein homotopy type.

\bigskip

This work was done during visits at IHES, MSRI, Weizmann Institute,
the Institut Fourier, ENS Paris and the Free University Berlin. The
author gratefully acknowledges the hospitality of these
institutions. The author is grateful to F.~Bogomolov, M.~Hedden,
C.~Livingstone and S.~Orevkov for helpful discussions, and to the
referee for pointing out inaccuracies and proposing a slight
simplification of the proof of Proposition 2. She would like to
thank C. Gourgues for typing a big part of the manuscript, M. Vergne
for drawing figures and Jo\"el Merker for correcting mistakes in
some of the figures.

\section{Examples}

The first two examples show that the complex curve bounded by the
link $L$ can be more complicated than the complex curve bounded by
the $\Omega '$-truncation $L'$
which occurs in part~1 of Theorem~2. The third example shows that
the statement of Theorem~1 does not extend to closures of
quasi-positive geometric braids $L$ that are $2$-braided links
around an analytic knot $K$. The example is due to S.~Orevkov
\cite{O1}.

\bigskip

\noindent {\bf Example 1.} (Twisted Whitehead doubles (winding number zero).)

\smallskip

Consider the analytic discs $\{ z_2 = 0 \} \cap {\mathbb B}^2$ and
$\{ z_1 = 0 \} \cap {\mathbb B}^2$. They intersect at the origin.
The union of their boundaries forms the Hopf link. Apply an
automorphism of the closed ball which maps the Hopf link to an
analytic link $L$ which is the union of two circles and is contained
in a 3-ball which is a subset of a small tubular neighbourhood
$N(K)$ of a given analytic knot $K$. Join a point $p$ on one circle
by a Legendrian arc in $N(K)$ with a point $q$ on the other circle.
Choose the arc without self-intersections and without intersection
points with the circles other than the endpoints. Moreover, the
Legendrian arc is chosen to be the longer part of a loop
representing a generator of the fundamental group of $N(K)$.
Consider a partition of the Legendrian arc into small closed arcs
with pairwise disjoint interior. For each small arc we take an
analytic disc on ${\mathbb B}^2$ (which extends to a complex curve
in a neighbourhood of the closed ball) such that its boundary lies
on the sphere and passes through the two endpoints of the arc and
through no other point of the large Legendrian arc. See figure~1.
Take an analytic function $f$ in a neighbourhood of
$\overline{{\mathbb B}^2}$ whose zero  set intersects the ball
$\mathbb{B}^2$ along the union of all these discs with the complex
curve bounded by the link $L$. For a suitable small complex number
$\varepsilon$ the set $\{ f = \varepsilon \} \cap {\mathbb B}^2$ is
a smooth complex curve with connected boundary. Part of the boundary
approximates a compact subset of $L \backslash ( \{ p \} \cup \{ q
\})$, the other part consists of two arcs close to the Legendrian
arc. The two arcs are traveled "{in opposite direction.}" This
follows from lemmas 3.6. and 3.7 in \cite{O2}.

\bigskip

$$
\includegraphics[width=5cm]{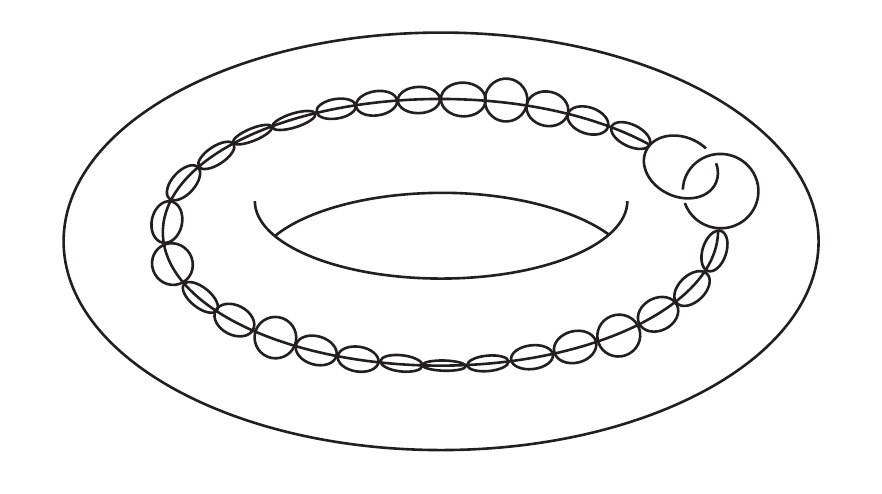}
$$

\centerline{Figure 1.}

\bigskip

\noindent {\bf Example 2.} (Sum of two analytic links.)

\smallskip

Take two analytic links $L_1$ and $L_2$ in the tubular neighbourhood
$N(K)$ of an analytic knot $K$. Join a point $p$ in one of the links
$L_1$ with a point $q$ in the other link $L_2$ by a Legendrian arc
in $N(K)$. The Legendrian arc is chosen without self-intersections
and with interior disjoint from the two links. As in example~1 we
find a complex curve $X$ in ${\mathbb B}^2$ which ``approximates''
the union of the complex curves bounded by the links and the
Legendrian arc. See figure~2. The boundary $\partial X$ is
connected. Part of it approximates a compact subset of $(L_1 \cup
L_2) \backslash (\{ p \} \cup \{ q \})$, the other part consists of
two arcs close to the Legendrian arc, the two arcs traveled in
opposite direction. One of the links can be taken to be an analytic
$n$-braided link around $K$, the second link may be obtained from an
arbitrary analytic knot in $\partial {\mathbb B}^2$ by an
automorphism of $\bar{\mathbb B}^2$ which maps the second link to a
$3$-ball contained in $N(K)$.

$$
\includegraphics[width=6cm]{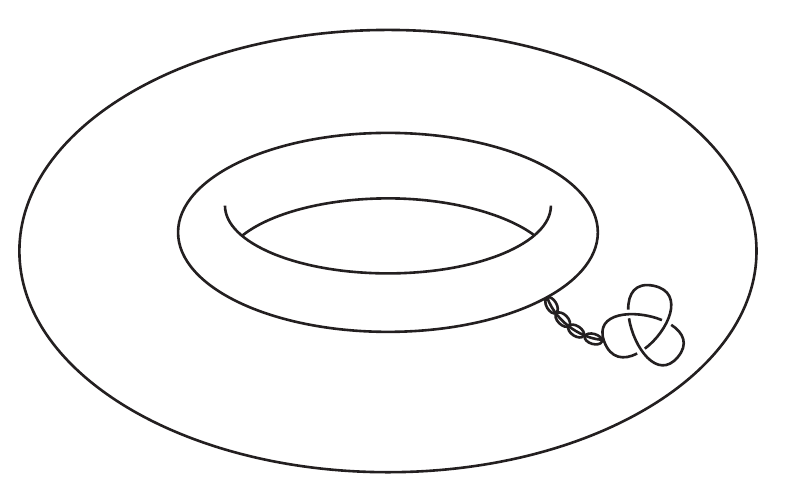}
$$

\centerline{Figure 2.}

\bigskip

\noindent {\bf Example 3.} (Closures of quasi-positive geometric
braids that are $2$-cables of closures of quasi-positive geometric
braids, see \cite{O1}, corollary~2.15.)

\smallskip
An $n$-braided link in a tubular neighbourhood of a knot $K$ is an
$n$-cable if it is isotopic in the tubular neighbourhood to an
$n$-braided link contained in the boundary of the tubular
neighbourhood. Consider the following braid $\Delta_n^2 \cdot
\sigma_{n-1}\cdot \ldots \cdot \sigma_1$ in the braid group
${\mathcal B}_n \,.$ Here $\sigma_1 , \ldots , \sigma_{n-1}$ are the
standard generators and $\Delta_n$ is Garside's half-twist. (By
induction $\Delta_0 = \Delta_1 = 1$, $\Delta_n = \sigma _1 \cdot
\ldots \cdot \sigma_{n-1} \, \cdot \Delta_{n-1}$.) The braid is
positive (i.e., a word containing only generators, not their
inverses), hence quasi-positive. The free isotopy class of its
closure is represented by a smoothly analytic knot $K_1 \subset
\partial {\mathbb B}^2$ (\cite{R}). (The closure of $\Delta_n^2$
is a link with $n$ connected components. Hence, the closure of the
considered braid is connected.) The (smooth) complex curve $X_1$
bounded by $K_1$ is (after adjusting near $\partial {\mathbb B}^2$)
a branched holomorphic covering of the disc with number of branch
points $b_1$ equal to the exponent sum of the braid (i.e., it is
equal to the sum of exponents of generators of ${\mathcal B}_n$
appearing in a representing word, see \cite{BO} or statement 2 of
proposition 4 below). Hence, $b_1$ is equal to $n \cdot (n-1) + n-1
= n^2 - 1$. The genus $g_4 (X_1)$ equals $\frac{b_1 - n+1}{2} =
\frac{n^2}{2} + O(n)$ by the Riemann-Hurwitz relation.

\smallskip

Orevkov proved \cite{O1}, Corollary~2.15, that for large $n$ and $N
\leq \frac{8}{3} \, n^2 + O(n)$ the braid $\sigma_1^{-N} \,
\Delta_{2n}^2 \in {\mathcal B}_{2n}$ is quasi-positive. Its closure
is a $2$-cable of $\Delta_n^2$. A modification of this example
provides a quasi-positive braid in ${\mathcal B}_{2n}$ whose closure
is a $2$-cable of $K_1$ and is connected. Indeed, consider
$\sigma_1^{-N} \, \cdot  
c(\sigma_{n-1}) \cdot \ldots \cdot c (\sigma_{1}) \, \cdot
\Delta_{2n}^2$ for odd $N =\frac{8}{3} \, n^2 + O(n)$. 
Here $c(\sigma_j) \overset{\rm def}{=} \sigma_{2j} \cdot
\sigma_{2j-1} \cdot \sigma_{2j+1} \cdot \sigma_{2j} \,$, $j =
1,\ldots , n-1$. See figure~3 for $n=3$. Denote the closure of this
braid by $K_2$. $K_2$ bounds a quasi-positive surface (see \cite{BO}
or statement 2 of proposition 4 below) with number of branch points
$b_2$ equal to the exponent sum of the braid, i.e., equal to $(2n)^2
- \frac{8}{3} \, n^2 + O(n) = \frac{4}{3} \, n^2 + O(n)$, and hence
of genus $\frac{2}{3} \, n^2 + O(n)$. Hence, for large $n$ the lower
bound for the $4$-ball genus is different from that in the case of
analytic satellites in small neighbourhoods of
smoothly analytic knots. 

$$
\includegraphics[width=7cm]{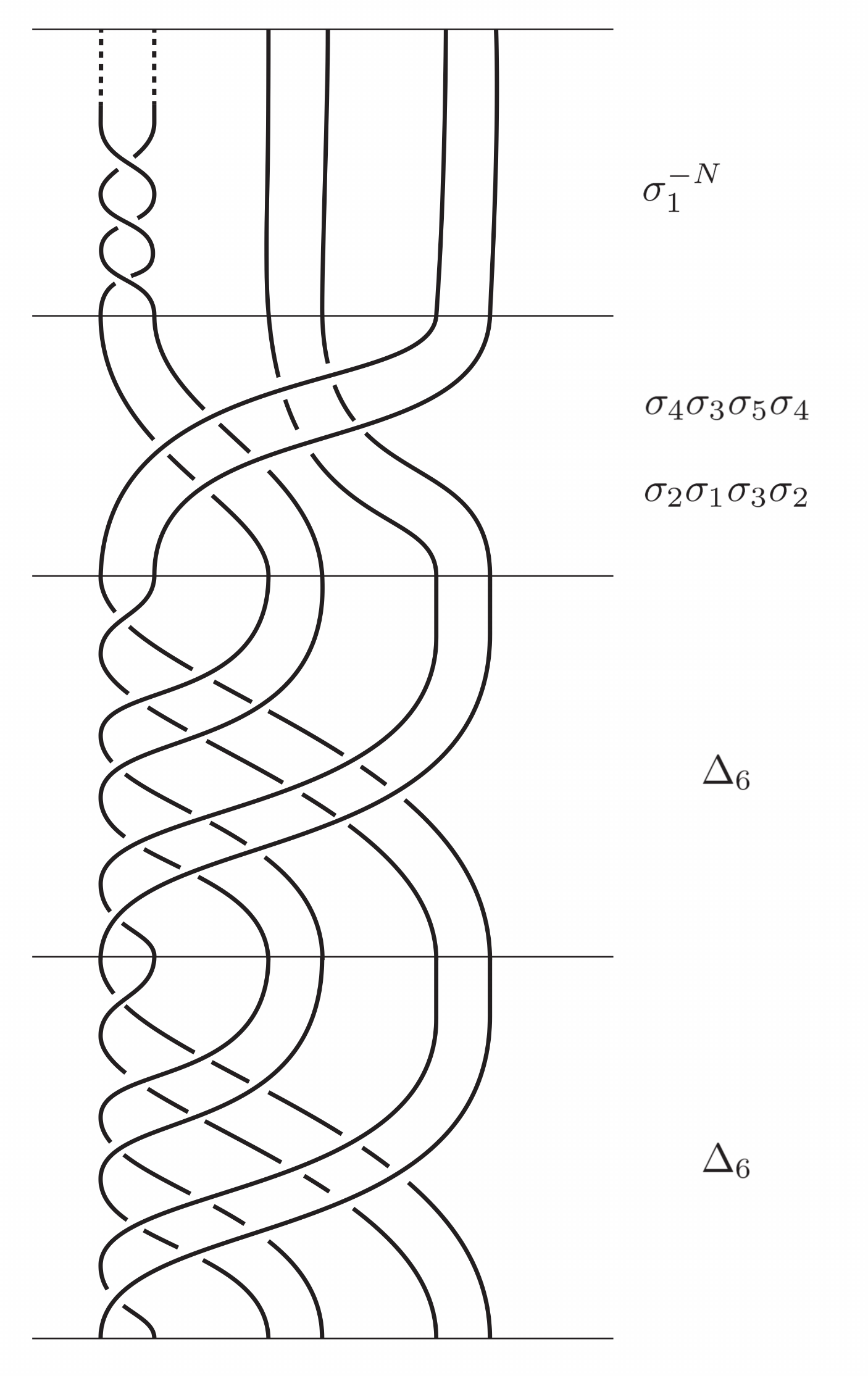}
$$

\centerline{Figure 3.}

\section{Proof of the $4$-ball genus estimates}

Theorem ~1 can be formulated in the more general situation when the
ball $\mathbb{B}^2$ is replaced by a relatively compact strictly
pseudoconvex domain $\Omega$ in a Stein surface $\tilde\Omega$ (for
short, when $\Omega$ is a Stein domain). In this section we will
work in this more general setting.

We start with the following simple but useful lemma.

\begin{lm} {\rm (\textbf{Tubular neighbourhood of knots and of complex
curves})} Let $\Omega$ be a Stein domain on a Stein surface
$\tilde\Omega$ and let $\tilde X$ be a relatively closed complex
curve in $\tilde\Omega$ (maybe, singular) which intersects $\partial
\Omega$ transversely along a knot $K$. Let $\tilde X \subset
\tilde\Omega$ be the zero set $\tilde X = \{ f =0 \}$ of an analytic
function $f$ in $\tilde\Omega$. For a relatively closed complex
curve $\tilde Y$ on $\tilde\Omega$ intersecting $\partial\Omega$
transversely along a link $L$, the inclusion $L \subset \{ \vert f
\vert < a \} \, \cap$ $\partial \Omega$ for some positive number $a$
implies the inclusion $Y = \tilde Y \cap \Omega \subset \{ \vert f
\vert < a \} \cap \Omega$.
\end{lm}

\bigskip

\noindent {\bf Proof of Lemma 2.} The inclusion $Y \subset \{ \vert
f \vert < a \} \cap \Omega$ follows from the maximum principle
applied to $f \mid Y$ and the fact that $\partial Y \subset \{ \vert
f \vert < a \} \cap \partial \Omega$. \hfill $\Box$

\bigskip

Suppose the gradient of $f$ does not vanish on $\tilde X \cap
\bar\Omega$. 
Then, possibly after taking for $\tilde\Omega$ a smaller Stein surface with $\Omega \Subset \tilde\Omega$ and choosing a small enough number $a > 0$, the gradient of $f$ does not vanish in a
neighbourhood of the subset $\{ \vert f \vert \le a \}$ of
$\tilde\Omega$ and there exists a holomorphic vector field $V^f =
(V_1^f , V_2^f)$ near this set such that $\frac{\partial}{\partial
z_1} \, f \cdot V_1^f + \frac{\partial}{\partial z_2} \, f \cdot
V_2^f = a$ (see e.g. \cite{H}, Theorem 7.2.9). For a domain
${\mathcal X} \Subset \tilde X$ with $\tilde X \cap \bar\Omega
\subset {\mathcal X}$ and for small enough $a$ the flow of $V^f$
defines a biholomorphic mapping $\phi^f$ from ${\mathcal X} \times
\, {\mathbb D}$ onto a tubular neighbourhood ${\mathcal T}_a
({\mathcal X})$ of ${\mathcal X}$. We obtain a trivial holomorphic
fiber bundle $ {\mathcal T}_a ({\mathcal X}) \to {\mathcal X}$ with
fiber $\phi^f (\{ p \} \times
 \, {\mathbb D})$ over the point $p \in {\mathcal X}$. We will
always consider the tubular neighbourhood ${\mathcal T}_a ({\mathcal
X})$ as total space of this bundle.

Let $a>0$ be small. Then there is a closed curve $K' \subset
{\mathcal X} \cap \Omega $ such that $K \cup - K'$ bounds an annulus
on ${\mathcal X}$ and for the domain ${X}' \subset {\mathcal X}$
bounded by $K'$ the tubular neighbourhood ${\mathcal T}_a
(\overline{X'})$ of the closure $\overline{X'}$ is contained in
$\Omega$. Here $-K'$ is obtained from $K'$ by inverting orientation.
The set $\mathcal{H}= {\mathcal T}_a (\partial{X'}) $ is a Levi-flat
hypersurface which is foliated into holomorphic discs $\Delta_z=
\phi^f(\{z\} \times \mathbb{D}),\, z \in
\partial{X'}$. It divides ${\mathcal T}_a ({\mathcal X})$
into two connected components. Denote by $A $ the connected
component which intersects $\partial \Omega$. For any a priory given
neighbourhood of $K$ in $\mathbb{C}^2$ the positive number $a$ can
be taken so small that the domain ${X}'$ can be chosen so that the
set $A$ is contained in the neighbourhood of $K$.

\begin{prop}
{\rm (\textbf{Truncation and closed braids})} Let in the described
situation $L$ be a link in $N(K)=\partial \Omega \cap \{|f| <a\}$
with winding number $w_{N(K)}(L)= n$. Suppose $L$ bounds a complex
curve $Y \subset \Omega$ and $\mathcal{H}$ intersects $Y$
generically. Then the following holds.

Either $n=0$ and then $ Y \cap \mathcal{H} = \emptyset$, or $n
> 0$. In the latter case  $L'  \overset {\rm def}= Y \cap \mathcal{H}$ is a closed $n$-braid in the solid torus
$\mathcal{H}$ around $K'$.
\end{prop}

\bigskip

\noindent {\bf Proof of Proposition 1.}  Since by Lemma 2 $ Y$ does
not meet the boundary $ \underset{z \in K'}{\bigcup} \partial
\Delta_z$ of $\mathcal{H}$ the intersection number of $ Y$ with
$\Delta_z$ is constant for $z \in K'$. For a neighbourhood $\tilde
X'$ of $\overline{{X}'}$ with ${\mathcal T}_a (\tilde X') \subset
\Omega$ the canonical projection $\rm{pr}:{\mathcal T}_a (\tilde X')
\to \tilde X' $ defines a (branched) holomorphic covering ${\rm pr}
\mid Y \cap {\mathcal T}_a (\tilde X')  \to \tilde X'$.

\smallskip

Since $\mathcal{H}$ and $Y$ intersect generically the latter
covering is unramified in a neighbourhood of the intersection
$\mathcal{H} \cap Y$. Orient $K'$ as boundary of ${X}'$ and orient
$L'$ as boundary of $Y \setminus \bar A$. Since a disc $\Delta_z$ is
either contained in $A$ or in $\partial A$ or does not meet
$\overline A$, the bundle projection ${\rm pr}$ maps
$\mathcal{T}_a(\tilde X' \cap A) $ into $\tilde X' \cap A$ and
${\mathcal T}_a (\tilde X')\setminus \bar A$ into $\tilde X'
\setminus \bar A$. Hence ${\rm pr} \mid Y \cap {\mathcal T}_a
(\tilde X')$ maps the side $Y \setminus \bar A$ of $L'$ on $Y$ to
the side $X \setminus \bar A$ of $K'$ on $X$, in other words ${\rm
pr} \mid L' : L' \to K'$ is orientation preserving if $L'$ and $K'$
are oriented as boundaries of complex curves in ${\mathcal T}_a
({\mathcal X}) \setminus \bar A$.

\smallskip

The intersection number of $Y$ with each disc $\Delta_z, \, z \in
\overline{{X}'}\,,\,$ equals $w_{N(K)} (L) = n$. This follows from
the fact that $L \cup -L'$ bounds the set $Y \cap A$ (after a small
perturbation of $Y$ we may assume that this set is a smooth
manifold). Indeed, for a neighbourhood of $\bar A$ which is
diffeomorphic to $K \times b^3$ let $e^{2\pi is}$, $s \in [0,1]$, be
the parameter in the direction of $K$. Integrate the form $ds$ along
$L$ and along $L'$ and apply Stokes' theorem. \hfill $\Box$

\medskip

The 4-ball genus estimates of Theorem 1 follow from the
Riemann-Hurwitz relation. Let $\textsf{X}$ and $\textsf{Y}$ be
smooth oriented surfaces or smooth oriented surfaces with boundary.
A smooth mapping ${ \textsf{p}} :  {\textsf{Y}} \to { \textsf{X}}$
is called a branched covering (opposed to holomorphic branched
covering) if it is an orientation preserving topological covering
outside the critical points, and has the form $\zeta \to \zeta ^n$
for a natural number $n \ge 2$ in suitable orientation preserving
complex coordinates on $\textsf{Y}$ in a neighbourhood of each
critical point. In case  $\textsf{X}$ and $\textsf{Y}$ have
non-empty boundary no critical point is allowed to be on the
boundary of $\textsf{Y}$. The Riemann-Hurwitz relation for branched
coverings of open Riemann surfaces is the following.

\bigskip

{\em Let $\textsf{X}$ and $\textsf{Y}$ be connected open Riemann
surfaces with smooth boundaries and let ${ \textsf{p}} : \bar
{\textsf{Y}} \to \bar{ \textsf{X}}$ be an 
$n$-fold branched covering 
Then}
\begin{equation}\label{eq3'}
\chi ({ \textsf{Y}}) = n \cdot \chi ({ \textsf{X}}) - \textsf{B} \,
.
\end{equation}
Here $\chi$ is the Euler characteristic and $\textsf{B}$ is the
number of branch points (counted with multiplicity).

Denote by $k({\textsf{X}})$ the number of boundary components of a
Riemann surface ${ \textsf{X}}$. Then
\begin{equation}\label{eq3''}
\chi ({ \textsf{X}}) = 2 - 2g ({ \textsf{X}}) - k ({ \textsf{X}}) \,
.
\end{equation}

The following proposition is a direct consequence of the
Riemann-Hurwitz relation and is used for the proof of Theorem 1.

\begin{prop} Let $\mathcal{X}$ be a connected open Riemann surface
with smooth connected boundary and let $\mathcal{Y}$ be an open
Riemann surface. 
Let $p: \overline{\mathcal{Y}} \to \overline{ \mathcal{X}}$ be an
orientation preserving $n$-fold branched covering.
Let $Y$ be a Riemann surface which contains $\mathcal{Y}$. Then

\begin{equation}\label{eq4'}
g(Y) \ge g(\mathcal{Y}) \ge n g(\mathcal{X}) - (n-1)
\end{equation}

\noindent The second inequality is an equality if and only if
$\mathcal{Y}$ is connected, the covering is unramified on
$\mathcal{Y}$ and $\mathcal{Y}$ has $n$ boundary components.

If $Y$ has connected boundary then

\begin{equation}\label{eq4''}
g(Y) \ge g(\mathcal{Y})
\ge  n g(\mathcal{X}) - [ \frac{n-1}{2}].
\end{equation}

\noindent If $n$ is odd then $g(Y) =  n g(\mathcal{X}) -
\frac{n-1}{2}$ if and only if 
$\,\mathcal{Y}$ and its boundary $\partial \mathcal{Y}$ are
connected, $g(Y)=g(\mathcal{Y})$ and the covering $\mathcal{Y} \to \mathcal{X}$ is
unramified.

\noindent If $n$ is even then the equality $g(\mathcal{Y}) =  n
g(\mathcal{X}) - [ \frac{n-1}{2}]$ can hold only if
the boundary $\partial \mathcal{Y}$ has one or two components.
\end{prop}

\noindent{\bf Proof.} Let  $m = m(\mathcal{Y})$ be the number of
connected components of $\mathcal{Y}$. Denote by ${\mathcal Y}_1 ,
\ldots , {\mathcal Y}_m$ the connected components of ${\mathcal Y}$.
Let $k_j$ be the number of boundary components of ${\mathcal Y}_j$,
let $n_j$ be the multiplicity of the covering $p \mid {\mathcal
Y}_j$ and let $B_j$ be the number of branch points (counted with
multiplicity) of the latter covering. Then $\mathcal{Y}$ has
$k(\mathcal{Y})= \sum_{j=1}^{m(\mathcal{Y})} \ k_j \geq 1\; $
boundary components , it has $\; B(\mathcal{Y}) =
\sum_{j=1}^{m(\mathcal{Y})}\ B_j
\geq 0 \;$ branch points, and has covering multiplicity $\,\; \sum_{j=1}^{m(\mathcal{Y})}\ n_j =n \; $. 
Use for each $j$ the relation (see \eqref{eq3'})
$$
\chi({\mathcal Y}_j)= n_j (1-2g({\mathcal X})) - B_j \, .
$$
Consider the sum over all $j$. 
We obtain
\begin{equation}\label{eq4'''}
\chi({\mathcal Y})= n (1-2g({\mathcal X})) - B({\mathcal Y}) \,.
\end{equation}
Apply \eqref{eq3''} to $\sf{X}=\mathcal{X}$, and to $\sf{Y}$ being a
connected component $\mathcal{Y}_j$ of $\mathcal{Y}$ with the number
of boundary components being $k_j$. Take the sum over all
$m(\mathcal{Y})$ connected components of $\mathcal{Y}$. We obtain
\begin{equation}\label{eq4'''a}
\chi(\mathcal{Y})= 2 m(\mathcal{Y}) -2 g(\mathcal{Y})
-k(\mathcal{Y}).
\end{equation}
Hence, since $k(\mathcal{Y}) \leq n$ and $m(\mathcal{Y}) \geq 1$ we
obtain from \eqref{eq4'''} and \eqref{eq4'''a}
\begin{equation}\label{eq4a}
2g({\mathcal Y}) = 
n(2g({\mathcal X}) -1) + B(\mathcal{Y}) +  2 m(\mathcal{Y}) - k(\mathcal{Y})\,     
\geq n \cdot 2g ({\mathcal X}) + 2(1-n) \,  ,
\end{equation}
and, hence, \eqref{eq4'}. The second inequality in \eqref{eq4'} is
an equality if and only if $B(\mathcal{Y})=0,\, m(\mathcal{Y})=1,\,$
and $k(\mathcal{Y})=n$.

Suppose $Y$ has connected boundary (hence, $Y$ is connected itself).
Then $Y \setminus \mathcal{Y}$ is connected and has
$k(\mathcal{Y})+1$ boundary components. Since
$\chi(Y)\,=\,\chi(\mathcal{Y})\, +\,\chi(Y \setminus \mathcal{Y})\,$
we obtain
$$
1-2g(Y) = \chi(Y) = \chi(\mathcal{Y}) +\chi(Y \setminus
\mathcal{Y})=
$$
\begin{equation}\label{eq4b}
n(1-2g(\mathcal{X}))-B(\mathcal{Y}) + 2-2g(Y \setminus
\mathcal{Y})-k(\mathcal{Y})-1.
\end{equation}
Hence,
\begin{equation}\label{eq4c}
g(Y) \ge ng(\mathcal{X})
-\frac{n}{2}+\frac{B(\mathcal{Y})}{2}+\frac{k(\mathcal{Y})}{2}.
\end{equation}
Since the right hand side of equation \eqref{eq4c} is an integral
number there is at least one branch point if $n-k(\mathcal{Y})$ is
odd. Hence
\begin{equation}\label{eq4d}
g(Y) \ge n \,g(\mathcal{X})- [\frac{n -k(\mathcal{Y})}{2}] \,.
\end{equation}
Since $k(\mathcal{Y})\ge 1$ we obtain \eqref{eq4''}.

It is clear (see \eqref{eq4b} ) that for odd $n$ equality
$g(Y)=n \,g(\mathcal{X})- \frac{n -1}{2}$ is attained if
and only if $k(\mathcal{Y})=1$ (i.e. $\partial \mathcal{Y}$ is
connected ), the covering is unramified (i.e.
$B(\mathcal{Y})=0$) and $g(Y)=g(\mathcal{Y})$. If $n$ is even the equality in \eqref{eq4} can
be attained only if 
$\,k(\mathcal{Y})=1 $ or $k(\mathcal{Y})=2\,$ (i.e. if $\,\partial
\mathcal{Y}\,$ has one or two connected components). The
proposition is proved. $\hfill \Box$

\medskip

\noindent {\bf Proof of statements 1 and 2 of Theorem 1 for smoothly
analytic knots $K$.} Proposition 1 implies
immediately that $n \geq 0$. 

Let 
$K = \tilde{ \mathcal{X}} \cap \Omega$ for a
smooth relatively closed complex curve $\tilde X$ in $\tilde
\Omega$, and let for a small number $a >0$  all values $\zeta$ with
$|\zeta| \leq a$ be regular values for the defining function $f$ of
$\tilde X$. Let $a$ be small so that for a relatively compact open
subset $X'$ of $X = \tilde X \cap \Omega$ with $X' $ diffeomorphic
to ${X}$ the inclusion $\mathcal{T}_a(\overline{X'}) \subset
\mathbb{B}^2$ holds.

Let $Y$ be the connected complex curve bounded by $L$. Then
$g_4(L)=g(Y)$. Indeed, let $Y_s$ be a smooth surface in the ball
$\mathbb{B}^2$ with the same boundary $\partial Y_s = \partial Y$ as
$Y$. If $Y_s$ is not connected we replace $Y_s$ by a connected
surface contained in the ball of the same genus with the same
boundary. This can be done by repeating the following procedure:
take a pair of connected components of $Y_s$, cut off a disc from
each of them, and glue back an annulus contained in $\mathbb{B}^2$
whose interior does not meet $Y_s$. Assume from the beginning that
$Y_s$ is connected.

We may assume after a small perturbation which does not change the
genus that $Y= \{F=0\} \cap \mathbb{B}^2$ for a polynomial $F$ in
$\mathbb{C}^2$ so that $Y$ extends to a smooth algebraic curve $C$ in
projective space $\mathbb{P}^2$. We repeat now the proof of
Corollary 1.3 of \cite{KM}. Let $D$ be a smooth curve of degree six
in $\mathbb{P}^2$ that meets $C$ transversally and is disjoint from
the ball. The double branched cover of $\mathbb{P}^2$ with branch
locus $D$ is a $K3$ surface. Since $C$ and $D$ intersect
transversally the inverse image $\hat C$ of $C$ under the covering
map is again a smooth complex curve.

Replace the lift of $Y$ to the $K3$ surface by the lift of $Y_s$ to
the $K3$ surface. We get a smooth connected surface $C_s$ in the
$K3$ surface in the same homology class as the curve $\hat C$. By
the theorem on the Thom conjecture the genus estimate
$$
g(\hat C) \leq g(C_s)
$$
holds. Equivalently, for the Euler characteristics we get the
estimate
$$
\chi(\hat C) \geq \chi(C_s)\,.
$$
Since the two surfaces coincide outside the lift of the ball
$\mathbb{B}^2$, we have
\begin{equation}\label{eq4e}
\chi(Y_s) \leq \chi(Y)\,.
\end{equation}

Since both, $Y$ and $Y_s$, are connected, the inequality
\eqref{eq4e} is equivalent to the inequality
\begin{equation}\label{eq4ee}
g(Y_s) \geq g(Y)\,.
\end{equation}

Inequality \eqref{eq4ee} shows that $g_4(L)=g(Y)$.

Apply Proposition 2 to  $Y$, to $X'$ instead of $\mathcal{X}$, and
to $Y \cap \mathcal{T}_a(X')$ instead of $\mathcal{Y}$.
\eqref{eq4ee} and \eqref{eq4''} give statement 1 of Theorem 1 for
smoothly analytic knots $K$, \eqref{eq4ee} and \eqref{eq4'} imply
statement 2 of Theorem 1 for smoothly analytic knots $K$. The case
when $K$ bounds a singular complex curve will be treated in the next
section. $\hfill \Box$

\bigskip

\section{Reduction to the case of $n$-braided links and holomorphic coverings}

Here we will prepare the proof of Theorem 2 and the proof of the
remaining statements of Theorem 1. 
For simplicity we assume that $\Omega$ and $\tilde{\Omega}$ are
strictly pseudoconvex balls in $\mathbb{C}^n$.

For 
the proof of Theorem 2 we will replace the Levi-flat hypersurface
$\mathcal{H}$ by a Levi-flat hypersurface which divides $\Omega$
rather than $\Omega \cap \{|f|<a\}$. The following lemma describes
the choice.

\begin{lm} {\bf (Deformation of tubular neighbourhoods of curves on strictly
pseudoconvex boundaries to Levi-flat hypersurfaces)} Let $\tilde X$
be a relatively closed complex curve in $\tilde\Omega$ (maybe,
singular) which intersects $\partial\Omega$ transversely along a
knot $K$. There are arbitrarily small tubular neighbourhoods
$N_{\partial} \subset
\partial\Omega$ of $K$ with the following properties. The boundary
$T = \partial N_{\partial}$ bounds a real-analytic Levi-flat
hypersurface $N_{\Omega} \subset \Omega$ which intersects $\tilde X$
transversely along a simple closed curve $K'$. Moreover, the
projection $N_{\Omega} \to K'$ defines a smooth fiber bundle with
fibers being analytic discs $\Delta_z,\; z \in K'$, with boundary
$\partial \Delta_z$ in $T \subset
\partial \Omega$. The union $N_{\partial} \cup T \cup N_{\Omega}$
bounds an open subset $A$ of $\Omega$ which is diffeomorphic to $S^1
\times b^3$ ($\,b^3$ a real $3$-dimensional ball), and the
pseudoconvex ball $\Omega'=\Omega \backslash \bar A$ is
diffeomorphic to $\Omega$. Moreover, $K \cup (-K')$ bounds an
annulus $\tilde X \cap A$ on $\tilde X$.

\end{lm}

$$
\includegraphics[width=6cm]{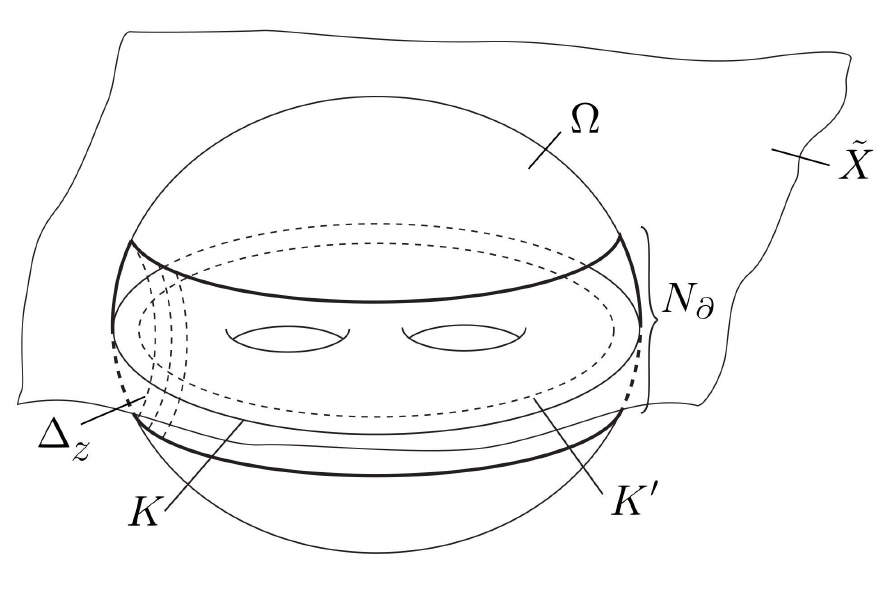}
$$

\centerline{Figure 4.}

\bigskip

Such tubular neighbourhoods $N_{\partial}$ will be called good.
Given any neighbourhood of $K$ in $\mathbb{C}^2$, it contains a good
neighbourhood $N_{\partial} \subset \partial \Omega$ of $K$ such
that the respective set $\bar A$ is contained in this neighbourhood.

\bigskip

\noindent {\bf Proof of Lemma 3.} Near each boundary point the
domain $\Omega$ is strictly convex in suitable holomorphic
coordinates, hence there exists a smooth fiber bundle over the curve
$K$ whose fiber over each $z \in K$ is a holomorphic disc embedded
into $\tilde\Omega$ which is complex tangent to $\partial\Omega$ at
the point $z$ and does not meet $\partial\Omega$ otherwise.
Moreover, we may assume that the union of the discs forms a smooth
Levi-flat hypersurface such that a neighbourhood of $K$ on
$\partial\Omega$ can be considered as the graph over this surface of
a function with non-degenerate quadratic form in the directions of
the holomorphic discs. By this fact we obtain a trivialization of
the smooth bundle over $K$ from Seifert framing on a tubular
neighbourhood of $K$ on $\partial \Omega$.

Note that the holomorphic discs are transversal to $\tilde X$.
Consider a $C^1$ approximation of the trivialized smooth bundle over
$K$ by a trivialized holomorphic disc bundle
$\mathcal{T}_{ct}(\mathcal{V}) \to \mathcal{V}$ over a neighbourhood
$\mathcal{V}$ of $K$ on $\tilde X$. We assume that $\mathcal{V}$ is
conformally equivalent to an annulus. 
If $\mathcal{V}$ is small enough the holomorphic disc through each
$z \in  {\mathcal{V}}' \overset{\rm{def}}= \mathcal{V} \cap \Omega$
intersects $\Omega$ along a connected simply connected set
$\Delta_z$. 
Choose a real analytic oriented loop $K'$ in ${\mathcal{V}}' =
\mathcal{V} \cap \Omega$ such that $K \cup (-K')$ bounds an annulus
on $\tilde X$. Put $N_{\Omega} = \underset{z \in K'}{\bigcup}
\Delta_z$, $T =
\partial N_{\Omega} \subset \partial \Omega$ and let $N_{\partial}$
be the connected component of $\partial\Omega \backslash T$
containing $K$. \hfill $\Box$

\bigskip

The following variant of Proposition 1 holds.

\smallskip

\noindent {\bf Proposition $1'$}. \textit{Let $\Omega$,
$\tilde\Omega$, $\tilde X$ and $K$ be as in Lemma 3. Let
$N_{\partial}$ be a good neighbourhood of $K$ and let $\Omega'$ be
the domain $\Omega \setminus \bar A$ for the set $A$ associated to
$N_{\partial}$ in Lemma 3. Then for any link $L$ contained in
$N_{\partial}$ which equals the transverse intersection of $\partial
\Omega$ with a relatively closed complex curve $\tilde Y$ in
$\tilde\Omega$ and has winding number $w_{N_{\partial}} (L) = n$ the
following holds.}

\smallskip

\textit{Either $n=0$ and then $\tilde Y \cap \partial \Omega' =
\emptyset$, or $n
> 0$. In the latter case (for a generic choice of $K'$ and, hence a generic choice of $\Omega'$)
the link $L' = \tilde Y \cap
\partial \Omega'$ is an $n$-braided link  (in $N_{\Omega} \subset
\partial \Omega'$) around $K' = \tilde X \cap
\partial \Omega'$.}

\textit{The statement also holds with $\Omega'$ replaced by a
strictly pseudoconvex domain $\Omega_1\,$  (depending on $K$ and
$L$), $\Omega' \subset \Omega_1 \subset \Omega,\, $ with $C^2$
boundary which is (away from corners of $\partial \Omega'$) $C^2$
close to $\partial{\Omega}'$.}

\medskip

The proof is a slight variation of the proof of Proposition 1 and is
left to the reader.

\bigskip

As a corollary we obtain a proof of statements 1 and 2 of Theorem 1
for $n=1$ in the case when $K$ bounds a singular complex curve.
Indeed, if $L$ is contained in a good neighbourhood $N_{\partial}$
of $K$ and $w_{N_{\partial}} (L)=1$ then $L'$ is isotopic to $K'$ in
$N_{\Omega} \subset \partial \Omega' $, hence $g(Y \cap \Omega') =
g_4(L')=g_4(K')=g_4(K)$ and $g_4(L)=g(Y) \ge g(Y \cap \Omega')$.

\bigskip

Suppose for some $a > 0\,$ the gradient of $f$ does not vanish on
the subset $\{ \vert f \vert < a \}$ of  $\,\tilde{\Omega}\,$, and
the set $\{ \vert f \vert < a \} \cap \partial\Omega$ is contained
in the good neighbourhood $N_{\partial} \subset \partial \Omega$ of
$K\,$ which appears in Lemma 3.

Let $\Omega'$ and $K'$ be as in Proposition $1'$. Put  $X' = \tilde
X \cap \Omega'$. Let $X''$ be a relatively compact open subset of
$X'$ such that $X' \setminus \overline{X''}$ is an annulus contained
in the set $\mathcal{V}'$ of the proof of Lemma 3. For the proof of
Theorem 2 we will need a smooth bundle $\mathcal{T}(\overline {X'})
\to \overline {X'}$
whose fibers are holomorphic discs. Moreover, the fibers over points
in a neighbourhood of $K'$ are contained in the discs $\Delta_z$ of
the bundle $\mathcal{T}_{ct}({\mathcal{V}}') \to {\mathcal{V}}'$,
and the total space of the bundle contains $\{|f|<a\} \cap
\overline{\Omega'}$.

For each relatively closed complex curve $\tilde Y$ in $\tilde
{\Omega}$ for which $\overline {Y'} \overset{\rm{def}}= \tilde Y
\cap \overline {\Omega'}$ is contained in $\{|f|<a\}$ the
intersection number of $\overline Y'$ with the discs of a bundle
with the listed properties
is constant. We will denote by $p=p_{\overline Y'}$ the restriction
of the bundle projection to $\overline Y'$.

\begin{prop} {\bf (Deformation of disc bundles)} There exists a smooth disc
bundle  $\mathcal{T}(\overline {X'}) \to \overline {X'}$  over
$\overline{X'}$ , which coincides over $X''$ with the bundle
$\mathcal{T}_{a''}(X'') \to X''$ for a small number $a''>0$, and
coincides over a small neighbourhood of $K'$ with a holomorphic disc
bundle whose fibers are subsets of the discs of the bundle
$\mathcal{T}_{ct}({\mathcal{V}}') \to {\mathcal{V}}'$. The fibers of
the bundle $\mathcal{T}(\overline {X'}) \to \overline {X'}$ are
holomorphic discs and the total space of the bundle contains
$\{|f|<a\} \cap \overline{\Omega'}$ for a small number $a>0$.

Moreover, for each relatively closed complex curve $\tilde Y$ in
$\tilde{\Omega}$ for which $Y= \tilde Y \cap \Omega$ is contained in
$\{|f|<a\}$, 
the induced mapping
$p_{\overline{Y'}}$ is a branched covering 
with branch locus outside a neighbourhood of $K'$.

Further, for a suitable strictly pseudoconvex domain $\Omega_1$, 
as in Proposition $1'$ (which depends on $\tilde X$ and $\tilde Y$) 
the bundle $\mathcal{T}(\overline {X'}) \to \overline {X'}$ can be
extended to a smooth disc bundle $\mathcal{T}(\overline{X_1}) \to
\overline{X_1}$ over $\overline{X_1}=\tilde X \cap
\overline{\Omega_1}$ such that the union of the fibers over $K_1=
\tilde{X} \cap
\partial \Omega _1$ constitute a neighbourhood of $K_1$ on
$\partial \Omega _1$ and the projection $\overline{Y_1} \setminus Y'
= \tilde Y \cap (\overline{\Omega_1} \setminus \Omega')  \to
\overline{X_1} \setminus X' $ induced by the bundle projection is an
unbranched orientation preserving covering. The total space
$\mathcal{T}(\overline X_1)$ of the bundle is contained in
$\overline {\Omega_1}$ and contains $\mathcal{T}_a(\tilde{X}) \cap
\overline{\Omega_1}$.


\end{prop}

We choose the trivialization of the bundle $\mathcal{T}(\overline
{X_1}) \to \overline {X_1}$ so that it induces Seifert framing on
$\mathcal{T}(K_1)$ for the knot $K_1= \tilde X \cap
\partial \Omega_1$ in $\partial \Omega_1$ .

\bigskip
\noindent {\bf Proof of Proposition 3.} 
Take a small number $a' >0$. For $z \in {\mathcal{V}}'$ the tangent
spaces at the point $z$ to the discs of both bundles ${\mathcal
T}_{a'} ({\mathcal{V}}') \to {\mathcal{V}}'$ and ${\mathcal T}_{ct}
({\mathcal{V}}') \to {\mathcal{V}}'$ are transversal to
${\mathcal{V}}'$. Hence, perhaps after shrinking the holomorphic
discs and using the previous notation for the bundles, we may
parametrize the discs so that the following holds. The disc of the
bundle ${\mathcal T}_{a'} ({\mathcal{V}}') \to {\mathcal{V}}'$
through $z \in {\mathcal{V}}'$ is parametrized by the natural
parametrization $\varphi^0_z(\zeta),\, \zeta \in \mathbb{D},$ as
integral curve of the vector field $V^f$. We parametrize the disc of the
bundle $\mathcal{T}_{ct}({\mathcal{V}}') \to {\mathcal{V}}' $
through $z \in {\mathcal{V}}'$ by $\varphi^1_z(\zeta),\, \zeta \in
\mathbb{D}$, so that for $j=1,2,,$ $\varphi^j_z(\zeta)= z + (v^j_z +
g^j_z(\zeta))\,\zeta\,,$ where $v^j_z$ are non-trivial vectors in
$\mathbb{C}^2$ which depend holomorphically on $z$ and
$g_z^j(\zeta)$ are holomorphic in $z$ and $\zeta$, and for all $z
\in \mathcal{V}'$ the Euclidean scalar product $(v^0_z,v^1_z)$
satisfies the inequality $|(v^0_z,v^1_z)| \ge c_0$ for a positive
constant $c_0$. Moreover, $|g_z^j| \le \frac{1}{2} c_0$. The
parametrization induces the previous trivialization of both bundles.
Indeed, the trivialization of the bundle ${\mathcal T}_{ct}
({\mathcal{V}}') \to {\mathcal{V}}'$  is related to Seifert framing
on a neighbourhood of $K$ on $\partial \Omega$ (see the proof of
Lemma 3). On the other hand the trivialisation of ${\mathcal T}_{a'}
({\mathcal X})$ is given by the normal vector field on ${\mathcal
X}$ that points in the direction of positive values of the function
$f$ and the set $\{ f
> 0 \} \cap \partial \Omega$ is, after a small generic
perturbation which fixes the part of the set which is contained in a
neighbourhood of $K$, a Seifert surface for $K$.

Take a small positive number $\varepsilon $. For each complex
parameter $t \in (1 + \varepsilon )\mathbb{D}$ the mapping
$G^t(z,\zeta)=\,t\,\varphi^0_z(\zeta) + (1-t)\, \varphi^1_z(\zeta) ,\,
z \in \mathcal{V}',\,\zeta \in \mathbb{D},\,$ is a holomorphic
diffeomorphism of $\mathcal{V}' \times \mathbb{D} $ onto a tubular
neighbourhood of $\mathcal{V}'$. The mappings $G^t$ and their
inverses $(G^t)^{-1}$ are uniformly bounded in the $C^1$ norm for $t
\in (1 + \varepsilon )\mathbb{D}$.

Take a smooth function $\,t(z), \, z \in \overline{X'} \setminus X''
\subset {\mathcal{V}}' \,,\,$ with values in
$\,(1+\varepsilon)\mathbb{D}\,,\,$ such that $t(z)=0$ in a
neighbourhood of $K'$ and $t(z)=1$ in a neighbourhood of $\partial
{X}''\,,\,$. Then for some small positive number $\delta$ the
holomorphic discs $t(z)\,\varphi^0_z(\zeta) + (1-t(z))\,
\varphi^1_z(\zeta) ,\, \zeta \in  \delta \mathbb{D},\,$ are the
fibers of a smooth disc bundle over $\overline{X'} \setminus X''$
with total space contained in $\overline{\Omega'}$. Over a small
neighbourhood of $\partial X''$ it coincides with the bundle
$\mathcal{T}_{a''}(X') \to X'$ with $a'' = \delta a'$. Over a
neighbourhood of $K'$ the fibers are subsets of the fibers of the
bundle $\mathcal{T}_{ct}({\mathcal{V}}') \to {\mathcal{V}}' $.
Extend the bundle to $X''$ by the holomorphic disc bundle
$\mathcal{T}_{a''}(X'') \to X''$. We obtain a smooth fiber bundle
$\mathcal{T}(\overline{X'}) \to \overline{X'}$ whose fibers are
holomorphic discs. Moreover, the total space of the obtained bundle
covers $\{|f|<a\} \cap \overline{\Omega'}$ for a small positive
number $a$.

Let $\tilde Y$ be a relatively closed complex curve in $\tilde
\Omega$ for which $Y=\tilde Y \cap  \Omega$ is contained in
$\{|f|<a\}$. Consider a transverse intersection point of
$\overline{Y'}$ with a disc of the bundle. Since the intersection is
positive the induced projection $p_{\overline{ Y'}}$ is an
orientation preserving diffeomorphism in a neighbourhood of the
intersection point. The points in a neighbourhood of $\partial Y'$
are transverse intersection points. Also, in a neighbourhood of
$\overline{X''}$ the bundle is holomorphic, hence $p_{\overline{
Y'}}$ is a branched holomorphic covering there.

Suppose $z_0$ is not in a neighbourhood of $K'$ or in a
neighbourhood of $\overline{X''}$ and  $\overline{Y'}$ intersects
the leaf which passes through the point $z_0$ of order $n \ge 2$. 
In coordinates $(z,\zeta)$ determined by the
diffeomorphism $G^{t(z_0)}$ the disc through $z_0$ equals the level
set $\{z=z_0\}$ and the equation of $Y'$ in a neighbourhood of the
intersection point is $\{z-z_0=a_n\,(\zeta -\zeta_0)^n +
\rm{higher\, order\, terms}\}$.

For $(z,t)$ close to $(z_0,t(z_0))$ the intersection of the disc
$\varphi_z^t(\zeta)=t \varphi^0_z(\zeta) + (1-t)
\varphi^1_z(\zeta)),\, \zeta \in \delta \mathbb{D},\,$ with $Y'$ is
described as follows. We have
\begin{equation}\label{eq5a}
\varphi_z^t(\zeta)= \varphi_z^{t(z_0)}(\zeta) + (t-t(z_0)) \,(
\varphi_z^0(\zeta)-\varphi_z^1(\zeta)),\, \zeta \in \delta
\mathbb{D}\,.
\end{equation}
The second term of the sum on the right equals $(t-t(z_0)) \,(
\varphi_z^0(\zeta)-\varphi_z^1(\zeta))= (t-t(z_0)) \zeta
g_z(\zeta),\,$ where $g_z(\zeta)$ is a holomorphic function in $z$
and $\zeta$ such that $|g_z(\zeta)| < c$ for a constant $c$ not
depending on the parameters. Apply the inverse $(G^{t(z_0)})^{-1}$
to \eqref{eq5a}. In the obtained coordinates $(z,\zeta)$ the disc
\eqref{eq5a} is given by
\begin{equation}\label{eq5b}
\zeta \to (z + (t-t(z_0)) \zeta h_1(z,\zeta, t-t(z_0)), \zeta +
(t-t(z_0)) \zeta h_2(z,\zeta, t-t(z_0))\,, \zeta \in \delta
\mathbb{D}\,,
\end{equation}
where $h_j,j=1,2,,\,$ are bounded holomorphic functions in
$(z,\zeta,\tilde t)$ with $z$ in a neighbourhood of $\overline{X'}
\setminus X''$, $\zeta \in \delta \mathbb{D}$, $\tilde t$ in a
neighbourhood of zero. Reparametrizing we may assume that 
$h_2$ is identically zero. For the intersection of $Y'$ with the
disc corresponding to $\varphi_z^t$ we obtain the equation
\begin{equation}\label{eq5c}
z-z_0 + (t-t(z_0)) \zeta h_1(z,\zeta, t-t(z_0))= a_n(\zeta -
\zeta_0)^n + \rm{higher\, order\, terms}\,.
\end{equation}
Let $t=t(z)$ for the chosen function $t(z)$ and suppose $\delta $ is
small and $z$ is close to $z_0$. Then for $|\zeta|< \delta$ the
point on the left of \eqref{eq5c} has distance at most
$\tilde{\delta} |z-z_0|$ from $z-z_0$ for a small constant
$\tilde{\delta}$. Hence for each $\zeta \in \delta \mathbb{D}$ the
function $z \to (z-z_0 + (t-t(z_0)) \zeta h_1(z,\zeta,
t-t(z_0)))^{\frac{1}{n}}$ is well-defined on the $n$-fold branched
cover of a neighbourhood of $z_0$ with branch locus $z_0$, and
$\zeta$ can be obtained as a quasiconformal function of
$(z-z_0)^{\frac{1}{n}}$. We proved that the induced mapping
$p_{\overline Y'}$ on $\overline Y'$ is a branched covering.

The last assertion of the proposition is clear. \hfill $\Box$

\bigskip

The following lemma is needed for the proofs of statement 4 of
Theorem 1 and of statement 3 of Theorem 2. It provides an isotopy of
an arbitrary smoothly analytic knot to a smoothly analytic knot with
more convenient properties.

\begin{lm}{\bf (Isotopy to closed geometric braids)} Let $K$ be a smoothly
analytic knot.  For small enough positive numbers $\epsilon$ there
exists an isotopy of $K$ to a smoothly analytic knot $K'=\mathcal{X}
\cap
\partial \mathbb{B}^2\,,$ where $\mathcal{X}$ is a relatively closed
curve in $(1+\epsilon)\mathbb{D} \times \epsilon \mathbb{D}$ such
that $\partial \mathcal{X}$ is a subset of $(1+\epsilon)\partial
\mathbb{D} \times \epsilon \mathbb{D}\,\,$ and $\mathcal{X} \cap \mathbb{B}^2$
is diffeomorphic to $\mathcal{X}$.  
\end{lm}

The lemma is a consequence of the following facts. Following
\cite{BO} a quasi-positive surface in $\mathbb{D} \times \mathbb{C}$
is a smooth proper embedding $\iota$ of a Riemann surface into
$\mathbb{D} \times \mathbb{C}$ with the following property.
For the canonical projection $P_1:\mathbb{D} \times \mathbb{C} \to \mathbb{D}$
the composition $P_1 \circ \iota$ is a branched covering.\\ 

Rudolph \cite{R} proved the following statement. \textit{A link $L$
in $\partial \mathbb{D} \times \mathbb{C}$ is the closure of a
quasi-positive geometric braid iff it bounds a quasipositive surface
in  $\partial \mathbb{D} \times \mathbb{C}.\,$ This happens iff $L$
is isotopic in $\partial \mathbb{D} \times
\mathbb{C}$ to the boundary of a complex curve.}\\

Bennequin \cite{B} proved that \textit{any oriented link in
$\partial \mathbb{B}^2$ which is positively transverse to the
complex tangent lines of $\partial \mathbb{B}^2$ is isotopic through
links with the same property (for short, it is transverse isotopic)
to a closed braid (in $\partial \mathbb{B}^2 \setminus \{z_1=0\}$).}

\bigskip

\noindent {\bf Sketch of proof of Lemma 4.} Let $X$ be the complex
curve in the ball bounded by $K$. Put $K_1=K$ and let $K_t$ be a
transverse isotopy in $\partial \mathbb{B}^2$ to a knot $K_2$ which
is contained in $\mathbb{C} \times \varepsilon \mathbb{D}$ for a
small positive number $\varepsilon$ and is the closure of a braid.
Consider the set
$$
X'=X \cup \bigcup_{t \in [1,2]} t K_t \cup \bigcup _{t \in
[2,\infty]} tK_2 \,.
$$

Here $tK_t \overset {\rm def}= \{ tz:\, z \in K_t\}$.  Boileau and
Orevkov \cite{BO} proved that, after smoothing, the set $X' \cap
\left(2 \mathbb{D} \times \mathbb{C}\right)$ is (diffeomorphic to) a
quasi-positive surface in $2\mathbb{D} \times \mathbb{C}$. By
Rudolph's theorem the link $L_0 \overset{\rm def}= X' \cap
\left(2\partial\mathbb{D} \times \mathbb{C}\right)$ is isotopic in
$2\partial\mathbb{D} \times \mathbb{C}$ to the boundary $L_1
\overset{\rm def}= \partial \mathcal{X}'$ of a complex curve
$\mathcal{X}'$. 
We may assume that $\mathcal{X}'$ is a subset of a quasipositive
surface in $(2+2\epsilon)\mathbb{D} \times \mathbb{C}$ such that the
projection to the first factor has branch locus in $2\mathbb{D}$. By
contraction in the $z_2$-direction the isotopy can be chosen so that
$L_1=\partial\mathcal{X}'$ (and hence also the complex curve
$\mathcal{X}'$) is contained in $\overline{2\mathbb{D}}
\times 2\epsilon \mathbb{D}$. 
The isotopy provides a family of links $L_t, \, t \in [0,1],\,$ in
$2\partial\mathbb{D} \times \mathbb{C}.\,$

Put $K'= \frac{1}{2}( \mathcal{X}' \cap 2\partial\mathbb{B}^2)$. It
remains to prove that $K_2$ is isotopic to $K'$ (since $K$ is
isotopic to $K_2$). Note that for a positive constant $\beta$ we
have the inclusion $K_2 \subset \{|z_1| > \beta\}$. Consider the
part $2
\partial \mathbb{B}^2 \cap \{|z_1|> 2\beta  \}$ of the sphere of
radius $2$. Assign to each point $z$ in this set the point of
intersection of the ray $\{t z:\, t >1\}$ with the set $\{ |z_1| =
2\} \times \mathbb{C}$. We obtain a diffeomorphism  $F$ of $2
\partial \mathbb{B}^2 \cap \{|z_1|> 2\beta  \}$ onto a subset of $\{
|z_1| = 2\} \times \mathbb{C}$ . Consider the
inverse $F^{-1}(L_t),\, t \in [0,1],\,$ of Rudolph's isotopy. 
This is an isotopy of links in $2\partial \mathbb{B}^2$ with
$F^{-1}(L_0)=2K_2$. If $\epsilon$ is small enough then $F^{-1}(L_1)
= F^{-1}(\partial\mathcal{X}') $ approximates $\mathcal{X}' \cap
2\partial \mathbb{B}^2$ well enough and, hence, is isotopic to it.
The lemma is proved. \hfill $\Box$


\medskip

Let $X$ be a smooth $2$-manifold (or a smooth $2$-manifold with
boundary). We call an embedding of a smooth $2$-manifold (or of a
smooth $2$-manifold with boundary) $Y$ into the product $X \times
\mathbb{D}$ an $n$-horizontal embedding if for the natural
projection $P_X: X \times \mathbb{D} \to X$ the restriction $P_X
\mid Y :Y \to X$ is an $n$-covering (smooth and unramified).

\begin{prop} Let $\epsilon $ be a small positive number.
Suppose an open Riemann surface $\mathcal{X}$ with smooth connected
boundary is holomorphically embedded into $(1+\epsilon) \mathbb{D}
\times \epsilon \mathbb{D}$ in such a way that the mapping $P_1 \mid
{\mathcal X} : {\mathcal X} \to (1+\varepsilon) \, {\mathbb D}$ is a
simple branched covering with branch locus in ${\mathbb D}$. Suppose
$i : \bar{\mathcal Y} \to \bar {\mathcal X} \times {\mathbb D}$ is a
smooth $n$-horizontal embedding of the closure of an
open Riemann surface ${\mathcal Y}$ into 
$\bar{\mathcal X} \times {\mathbb D}$. Then the following statements hold.\\
\begin{enumerate}
\item[1.]{\bf (Isotopy of $n$-horizontal embeddings to holomorphic
embeddings)} There exists a simply connected smoothly bounded domain
${\mathcal D} \subset {\mathbb D}$ containing the branch locus of
$P_1 \mid {\mathcal X}$ such that the following holds.  Put $X =
(P_1 \mid {\mathcal X})^{-1} ({\mathcal D}) = {\mathcal X} \cap
({\mathcal D} \times {\mathbb C})$ and $Y = 
\mathcal{Y} \cap (X \times \mathbb{D})\,.\,$ The embedding of the
two-manifold $Y$
into $X \times \mathbb{D}$ is isotopic through horizotal embeddings
to a holomorphic embedding of a Riemann surface into $X \times
\mathbb{D}$.

\item[2.] {\bf (Isotopy classes of boundary links)} 
Denote by $w_{\mathcal Y} \in {\mathcal B}_n$ a braid whose
conjugacy class corresponds to the isotopy class of the boundary
link $\partial {\mathcal Y} \subset
\partial {\mathcal X} \times {\mathbb D}$ of statement 1. Let $w \in {\mathcal
B}_n$ be a quasi-positive braid. Then with $\mathcal{D}$ as in
statement 1 there exists a domain $\mathcal{D}_1,\; \mathcal{D}
\subset \mathcal{D}_1 \subset
\mathbb{D},\,$ such that with $X_1= (P_1)^{-1}(\mathcal{D}_1)$ 
there exists a smooth embedding $i : \bar Y_1 \to \bar X_1 \times
{\mathbb D}$ of the closure $\bar Y_1$ of an open Riemann surface
$Y_1$  which is holomorphic on the Riemann
surface $Y_1$ and such that $P_{\mathcal X} \mid Y_1$ is a branched
$n$-covering of $X_1$ and the isotopy class of the boundary link
$\partial Y_1 \subset
\partial X_1 \times \mathbb D$ corresponds to the conjugacy class of
$w \cdot w_{\mathcal Y}$. The number of branch points equals the
exponent sum of the braid $w$.
\end{enumerate}
\end{prop}

Notice that the complex curves $X$ and $\mathcal{X}$ are
diffeomorphic. By a contraction in the $z_2$-direction we may also
assume that $\mathcal{X}$ and $\mathcal{X} \cap \mathbb{B}^2$ are
diffeomorphic.

Postpone the proof of Proposition ~4.

The following lemma relates the boundaries of the embedded surfaces
in Proposition ~4 to analytic links which are braided links around
analytic knots.
\bigskip

\begin{lm}\textbf{(Boundary links of horizontally embedded surfaces and braided links around knots)}
Let $\epsilon>0$ be a sufficiently small number, and let ${\mathcal
X}$ be a Riemann surface which is holomorphically embedded into
$(1+\epsilon)\mathbb{D} \times \epsilon \mathbb{D}$ so that $P_1
\mid {\mathcal X} : {\mathcal X} \to (1+ \varepsilon) \, {\mathbb
D}$ is a branched covering with branch locus contained in ${\mathbb
D}$ such that $\mathcal{X}$ is diffeomorphic to $\mathcal{X} \cap
\mathbb{B}^2$. Denote by $K$ the knot $K = {\mathcal X} \cap
\partial {\mathbb B}^2$.

\smallskip

Let ${\mathcal D} \subset {\mathbb D}$ be a simply connected
smoothly bounded domain containing the branch locus of $P_1 \mid
\mathcal{X}$ so that $X = {\mathcal X} \cap ({\mathcal D} \times
{\mathbb C})$ is diffeomorphic to ${\mathcal X}$. Suppose there is
an open Riemann surface $Y$ and a holomorphic embedding $i : Y \to
{\mathcal T}_a (X)$ into a small tubular neighbourhood of $X$ such
that for the projection $P_X:{\mathcal T}_a (X)  \to X $ the mapping
$p=P_X \circ i$ is a branched holomorphic $n$-covering.

\smallskip

Then there exists an isotopy of $K$ in $\partial \mathbb{B}^2$ to a
smoothly analytic knot $\tilde K$ such that any a priori given
tubular neighbourhood of $\tilde K$ contains an analytic link
$\tilde L$ which is an $n$-braided link around $\tilde K$ with
pattern equal to the isotopy class of $\partial Y $ in $ {\mathcal
T}_a(\partial X)$. The complex curve $\tilde Y$ in $\mathbb{B}^2$
bounded by $\tilde L$ is diffeomorphic to $Y$.
\end{lm}

\noindent {\bf Proof of Lemma 5.} Let ${\mathcal D}_t$, $t = [0,1],$
be a continuous decreasing family of simply connected smoothly
bounded domains with ${\mathcal D}_0 = {\mathbb D}$ and ${\mathcal
D}_1 = {\mathcal D}$. Let $\alpha$ be a continuous function on
$[0,1]$ with $\alpha(0)=1$ and $\alpha(t) > 1$ for $t>0$. Denote by
$\omega_t : {\mathcal D}_t \to \alpha(t) \,{\mathbb D}$ the
conformal mapping with $\omega'_t (0) > 0$. The $\omega_t$ depend
continuously on $t$.

\smallskip

Let $t \to s(t),\; t \in [0,1],\,$ be a continuous decreasing
positive function with $s(0)=1$ which takes sufficiently small
values for $t$ away from $0$. For $z=(z_1,z_2) \in \mathcal{D}_t
\times \mathbb C, \; t \in [0,1],\;$ we put
$\mathcal{G}_t(z)=(\omega_t(z_1),\, s(t)\cdot z_2)$. The mapping
$\mathcal{G}_t \mid \mathcal{X} \cap (\mathcal{D}_t \times \mathbb
C)\,$ is a conformal map onto a (relatively closed) complex curve
$X_t$ in $\alpha(t)\,{\mathbb D} \times \epsilon{\mathbb D}$. If
$\epsilon$ is small enough then for suitable choices of the
functions $\alpha$ and $s$ the intersection of each $X_t$ with
$\partial {\mathbb B}^2$ is transversal, and, hence, $K_t = X_t \cap
\partial {\mathbb B}^2$, $t \in [0,1]$, is a transversal isotopy.
Moreover, if $s(1)>0$ is small then  the subset $\tilde X = X_1 \cap
{\mathbb B}^2$ of $X_1$ is close to $X_1$ and is diffeomorphic to
$X_1$. Put $\tilde K =
\partial \tilde X$.

\smallskip

We may assume that $a > 0$ is sufficiently small and $Y$ is
identified with an embedded submanifold of ${\mathcal T}_a (X)$. The
Riemann surface $Y_1=\mathcal{G}_1(Y)$ is conformally equivalent to
$Y$ and embedded into the tubular neighbourhood
$\mathcal{G}_1({\mathcal T}_a (X))$ of $X_1$. For a suitable choice
of the function $\alpha$ and small $a$ the Riemann surface $\tilde Y
= Y_1 \cap {\mathbb B}^2$ is diffeomorphic to $Y_1$ and the boundary
$\tilde L \overset{\rm def}{=}\partial \tilde Y \subset \partial
{\mathbb B}^2$ is an $n$-braided link around $\partial {\tilde X}$
contained in the small tubular neighbourhood
$\mathcal{G}_1({\mathcal T}_a (X)) \cap
\partial {\mathbb B}^2$ of $\tilde K$. Moreover, the isotopy class of
$\partial \tilde Y $ in $\mathcal{G}_1({\mathcal T}_a (X)))$ equals
the pattern of $\tilde L$.  \hfill $\Box$

\bigskip

For the proof of Proposition 4 it will be convenient to work with
the following terminology.

A continuous mapping from a smooth manifold (or a smooth manifold
with boundary) $X$ into the set of all monic polynomials $\overline
{\mathfrak{P}_n}$ of degree $n$ is a quasi-polynomial of degree $n$.
It can be written as function in two variables $x \in X,\; \zeta \in
\mathbb{C}, \;$ i.e. $\,\mathcal{P}(x,\zeta)= a_0(x) + a_1(x)\zeta+
... + a_{n-1}(x)\zeta^{n-1}
+ \zeta^n, \;$
for continuous functions $a_j, \; j=1,...,n\,,$ on $X$. If the image
of the map is contained in the space $\mathfrak{P}_n$ of monic
polynomials of degree $n$ without multiple zeros, it is called
separable. 
For a separable
quasipolynomial, considered as a function $\mathcal{P}$ on $\,X
\times \mathbb{C}\,$, its zero set $\,\mathfrak{S}_{\mathcal{P}}=
\{(x,\zeta) \in X\times \mathbb{C},\; \mathcal{P}(x,\zeta)=0\}\;$ is
a surface which is $n$-horizontally embedded into $X \times
\mathbb{C}^n$. Vice versa, each $n$-horizontally embedded
$2$-manifold in $X\times \mathbb{C}$ corresponds to a separable
quasi-polynomial of degree $n$ on $X$. An isotopy of separable
quasi-polynomials (i.e. a family of separable quasi-polynomials
depending continuously on a parameter in $[0,1]$) is equivalent to
an isotopy of $n$-horizontally embedded manifolds.

Notice that the set $\mathfrak{P}_n$ is biholomorphic to
$\mathbb{C}^n \setminus \{\sf{D}_n =0 \},\,$ where $\sf{D}_n$
denotes the discriminant.

\bigskip

\noindent {\bf Proof of Proposition 4.} For the {\bf proof of
assertion 1} we consider a simple smooth arc $\Gamma$ (i.e. a
diffeomorphic image of the closed unit interval) in the disc
$\mathbb{D}$ which passes once through each point of the branch
locus $E$ of the covering $\mathcal{X} \to (1+\epsilon) \mathbb{D}$.
Let $\Gamma_{\mathcal{X}}
= \mathcal{X} \cap (P_1)^{-1}(\Gamma) \subset \mathcal{X}$ 
and let $\mathcal{P}$ be the mapping from $\bar {\mathcal{X}}$ into
$C_n(\mathbb{C}) \diagup \mathcal{S}_n \cong \mathfrak{P}_n \cong \mathbb{C}^n \setminus
\{\sf{D}_n =0 \} \subset \mathbb{C}^n\,$ defined by the embedding of
$\bar{\mathcal{Y}}$ into $\bar {\mathcal{X}} \times \mathbb{D} $.
The subset $\Gamma_{\mathcal{X}}$ of $\mathcal{X}$ has no interior
point. Hence, by Mergelyan's Theorem for Riemann surfaces \cite{K}
the restriction $\mathcal{\mathcal{P}} \mid \Gamma_{\mathcal{X}}$
can be approximated uniformly on $\Gamma_{\mathcal{X}}$ by an
analytic mapping of a neighbourhood $X \subset \mathcal{X}$ of
$\Gamma_{\mathcal{X}}$ into $\mathbb{C}^n$. After perhaps shrinking
$X$, we may assume that the restriction $P_1 \mid X$ is a branched
covering onto a domain $\mathcal{D} \subset \mathbb{D}$. The set
$\Gamma_{\mathcal{X}}$ is a deformation retract of $\mathcal{X}$.
After shrinking $X$ we may assume that $X$ is
diffeomorphic to $\mathcal{X}$. 

If the approximation is good enough then, after, perhaps, shrinking
$X$ again, the image of $X$ under the approximating mapping is
contained in the symmetrized configuration space. We obtain a
holomorphic mapping of $X$ into the symmetrized configuration space
which is isotopic to $\mathcal{P} \mid X$ through smooth mappings
into symmetrized configuration space. The isotopy of mappings
defines an isotopy of $n$-horizontal embeddings into $X \times
\mathbb{D}$ which joins the embedding of $\mathcal{Y} \cap \left(X
\times \mathbb{D}\right)$ with a holomorphic $n$-horizontal
embedding of a complex curve into $X \times \mathbb{D}$. We proved
assertion 1.

For the {\bf proof of assertion 2} we need Lemma~6 below. The
following construction prepares its statement. Let ${\mathcal X}$ be
an open Riemann surface with smooth connected boundary. Take any
smoothly bounded domain ${\mathcal X}_0 \subset {\mathcal X}$ which
is a strong  deformation retract of ${\mathcal X}$, and take simply
connected smoothly bounded domains ${\mathcal X}_j \Subset {\mathcal
X} \backslash {\mathcal X}_0$, $j = 1,\ldots , k$, with pairwise
disjoint closure. Consider simple smooth pairwise disjoint arcs
$\gamma_j : [0,1] \to {\mathcal X},\,\,j=1,\ldots , k,\,$ such that
for each $j$ the interior of $\gamma_j$ is contained in ${\mathcal
X} \backslash \underset{j=0}{\overset{k}{\bigcup}} \bar{\mathcal
X}_j$ and $\gamma_j$ joins a boundary point of ${\mathcal X}_0$ with a
boundary point of ${\mathcal X}_j$. Consider disjoint "rectangles"
$R_j$ around the $\gamma_j$, {\em i.e.} for some $\varepsilon > 0$
the $R_j$ are diffeomorphic mappings of $(-\varepsilon ,
\varepsilon) \times [0,1]$ with image in ${\mathcal X}$, 
and such that $R_j \mid \{ 0 \} \times [0,1]$ equals
$\gamma_j$ and the sides $R_j ((-\varepsilon , \varepsilon) \times
\{ 0 \})$ and $R_j ((-\varepsilon , \varepsilon) \times \{ 1 \})$
are contained in $\partial {\mathcal X}_0$, respectively in
$\partial {\mathcal X}_j$. Suppose the domain ${\mathcal X}_0 \cup
\underset{j=1}{\overset{k}{\bigcup}} {\mathcal X}_k \cup
\underset{j=1}{\overset{k}{\bigcup}} R_j ((-\varepsilon ,
\varepsilon) \times [0,1])$ has smooth boundary (see fig.~5). Then
it is again a deformation retract of ${\mathcal X}$. We call any
domain of the above described type a thickening of
$\underset{j=0}{\overset{k}{\bigcup}} {\mathcal X}_k \cup
\underset{j=1}{\overset{k}{\bigcup}} \gamma_j$.

\bigskip

$$
\includegraphics[width=8cm]{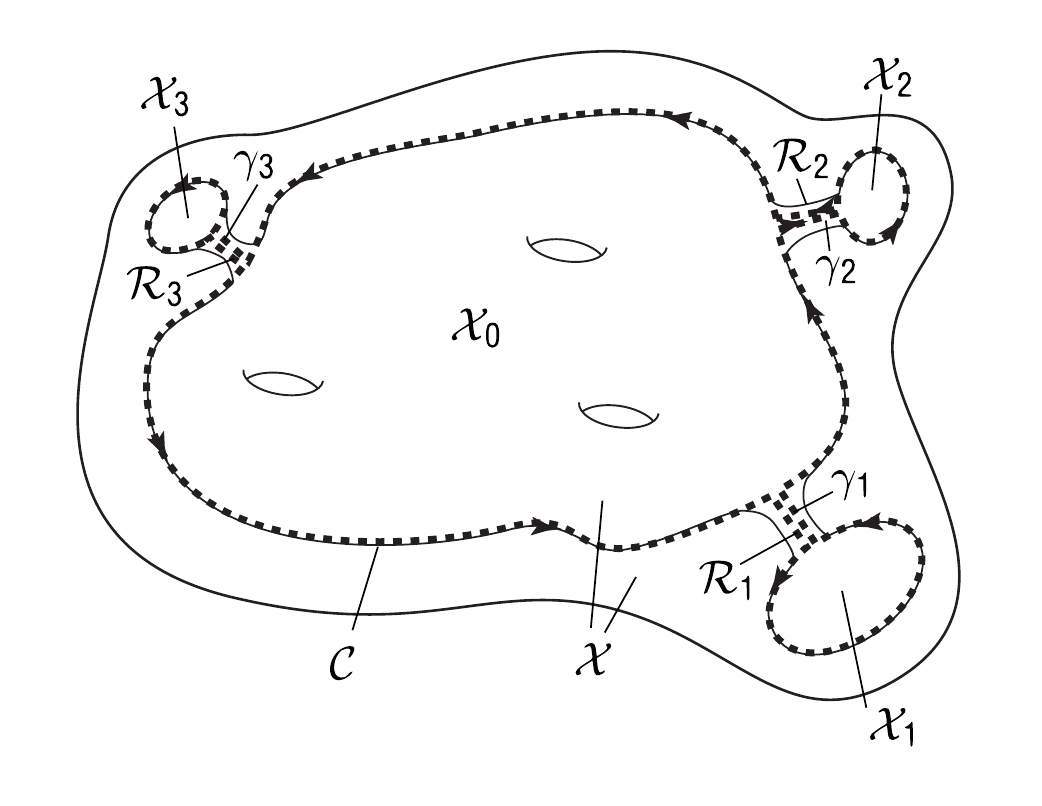}
$$

\centerline{Figure 5.}

\bigskip

Denote by $P_{\mathcal X} : {\mathcal X} \times {\mathbb D} \to
{\mathcal X}$ the canonical projection 
to the first factor. Consider a branched holomorphic covering $p$
from an open Riemann surface $\mathcal{Y}$ onto $\mathcal{X}$, and a
holomorphic embedding $i$ of $\mathcal{Y}$ into the disc bundle
$\mathcal{X} \times \mathbb D$ such that $P_{\mathcal X} \circ i =
p$ (for short $i$ lifts $p$). The following Lemma~6 provides the
embedding of another open Riemann surface $Y$ into the holomorphic
disc bundle over a deformation retract $X$ of $\mathcal{X}$ such
that the embedding lifts a branched covering with "more" branch
points and has the following property: the isotopy class of the
embedding of the boundary $\partial Y$ into $\partial X \times
\mathbb D$ differs from that of $\partial \mathcal{Y} \subset
\partial X \times \mathbb {D }$ by a prescribed quasi-positive braid.

\begin{lm}{\bf (Adding branch points)}
Let $\mathcal{X}$ and $\mathcal{Y}$ be open Riemann surfaces with
smooth boundary. Suppose the boundary of  $\mathcal{X}$ is
connected. Let $i : \overline{\mathcal Y} \to \overline{\mathcal X}
\times {\mathbb D}$ be a smooth embedding which is holomorphic on
${\mathcal Y} $ and such that for the canonical projection
$P_{\overline{\mathcal X}} :\overline {\mathcal X} \times {\mathbb
D} \to \overline{\mathcal X}$ the restriction $P_{\overline{\mathcal
X}} \mid {\mathcal Y}$ is a ( simple branched or unbranched)
holomorphic $n$-covering of ${\mathcal X}$. Denote by $w_{\mathcal
Y} \in {\mathcal B}_n$ a braid whose conjugacy class corresponds to
the isotopy class of the boundary link $\partial {\mathcal Y}
\subset
\partial {\mathcal X} \times {\mathbb D}$. Let $w \in {\mathcal
B}_n$ be a quasi-positive braid of exponent sum $m$. Let ${\mathcal
X}_0$ be as above. Suppose ${\mathcal X} \backslash {\mathcal X}_0$
does not intersect the branch locus of ${\mathcal P}_{\mathcal X}
\mid {\mathcal Y}$.

Then for $k \ge m$ and for each collection ${\mathcal X}_0 ,
{\mathcal X}_1 , \ldots , {\mathcal X}_k, \gamma_1 , \ldots ,
\gamma_k,$ as above there exists a thickening $X \subset {\mathcal
X}$ of $\underset{j=0}{\overset{k}{\bigcup}} {\mathcal X}_j \cup
\underset{j=1}{\overset{k}{\bigcup}} \gamma_j$ and a holomorphic
embedding of an open Riemann surface $i : Y \to X \times {\mathbb
D}$ into the disc bundle over $X$ such that $P_{\mathcal X} \mid Y$
is a simple branched $n$-covering of $X$ whose branch locus contains
exactly one point in each member of a chosen collection of $m$ sets
among the $\mathcal{X}_j,\,j=1,\ldots,k,\;$ and may contain also
some points in $\mathcal{X}_0$.  Moreover, the isotopy class of the
boundary link $\partial Y \subset
\partial X \times {\mathbb D}$ corresponds to the conjugacy class of the braid
$w \cdot w_{\mathcal Y} \in {\mathcal B}_n$. The same conclusion
holds for any thickening of $\underset{j=0}{\overset{k}{\bigcup}}
{\mathcal X}_j \cup \underset{j=1}{\overset{k}{\bigcup}} \gamma_j$
which is contained in $X$.
\end{lm}

\noindent {\bf End of Proof of assertion 2 of Proposition 4.} Let
$\mathcal{D}$ be the domain chosen in assertion 1, let $d_j,\,
j=1,\ldots, m,\,$ be pairwise disjoint smoothly bounded open
topological discs contained in $\mathbb{D} \setminus
\mathcal{\mathcal{D}},\, $ and let $\Gamma_j,\, j=1,\ldots,n,\,$ be
pairwise disjoint arcs with interior contained in $\mathbb{D}
\setminus (\bar{\mathcal{D}} \cup \bigcup \bar {d_j})$ such that
$\Gamma_j$ joins a boundary point of ${\mathcal{D}}$ with a boundary
point of $d_j$. Let $\mathcal{D}_1 \subset \mathbb{D}$ be a small
enough thickening of $\mathcal{D} \cup \bigcup (d_j \cup \Gamma_j)$.
Then $(P_1)^{-1}(\mathcal{D}_1)$ is a small thickening of
$(P_1)^{-1}(\mathcal{D} \cup \bigcup (d_j \cup \Gamma_j))\,.$  Apply
Lemma 6, where we take $\mathcal{X}_0= (P_1)^{-1}(\mathcal{D})
\subset \mathcal{X}$,  $\mathcal{X}_j, \,\, j=1,\ldots,k,\,$ running
over all  components of the preimages under $P_1$ of the $d_i$ for
all $i$, and $\gamma_j ,\,\, j=1,\ldots,k,\, $ (with the respective
label) running over all components of the preimages of all
$\Gamma_i$. Choose the collection of $m$ of the $\mathcal{X}_j$'s so
that $P_1$ is injective on the union of the sets of the collection.

Lemma 6 implies statement 2 of Proposition 4. \hfill $\Box$

\bigskip

\noindent {\bf Proof of Lemma 6.} Denote by $\mathcal{A}$ the
compact set $\mathcal{A} = \underset{j=0}{\overset{k}{\bigcup}}
\bar{\mathcal X}_j \cup \underset{j=1}{\overset{k}{\bigcup}}
\gamma_j$. Let ${\mathcal C}$ be the (oriented) closed curve which
is contained in the boundary $\partial \mathcal{A}$ and surrounds
$\mathcal{A}$ counterclockwise (see fig.~5). ${\mathcal C}$ is
obtained from $\partial {\mathcal X}_0$ (traveled counterclockwise)
in the following way: for each $j=1,\ldots , k,$ we cut $\partial
\mathcal{X}_0$ at the starting point $q_j$ of $\gamma_j$ and insert
the curve $\gamma_j$ followed by $\partial {\mathcal X}_j$ traveled
from the endpoint of $\gamma_j$ surrounding ${\mathcal X}_j$
counterclockwise, followed by $-\gamma_j$ ($\gamma_j$ with inverse
orientation).

\smallskip

Consider the closed geometric braids $\partial\mathcal{Y} \subset
\partial \mathcal{X} \times \mathbb{C}$ and $\mathcal{Y} \cap
\left(\partial \mathcal{X}_0 \times \mathbb{\mathbb{C}}\right).\, $
Let $Q$ be a point in $\partial\mathcal{X}$ and let $q =q_1 \in
\partial \mathcal{X}_0$ be the initial point of the arc
$\gamma_1.\,$ Put $E_n^Q= P_{\mathbb{C}}\left(\mathcal{Y}\cap(\{Q\}
\times \mathbb{C})\right)$ and $E_n^q=
P_{\mathbb{C}}\left(\mathcal{Y}\cap(\{q\} \times
\mathbb{C})\right)$. Here $P_{\mathbb{C}}:\mathcal{X} \times
\mathbb{C} \to \mathbb{C}$ is the canonical projection.

Denote by $\Gamma$ an arc 
whose interior parametrizes $\partial\mathcal{X} \setminus \{Q\}$
(with orientation induced from the orientation of
$\partial\mathcal{X}$ as boundary of $\mathcal{X}$) and whose
endpoints are equal to $Q$. Let $\Gamma_0$ be the respective object
for $\partial\mathcal{X}_0$ and the point $q$. The closed geometric
braid $\partial\mathcal{Y}$ defines a continuous map from $\Gamma$
to $\mathfrak{P}_n$ with base point $E_n^Q$, the closed geometric
braid $\mathcal{Y} \cap (\partial\mathcal{X}_0 \times \mathbb{C})$
defines a continuous map from $\Gamma_0$ to $\mathfrak{P}_n$ with
base point $E_n^q$.

Let $\gamma$ be a closed arc with interior contained in ${\mathcal
X} \backslash \left(\underset{j=0}{\overset{k}{\bigcup}}
\bar{\mathcal X}_j \cup \underset{j=1}{\overset{k}{\bigcup}}\gamma_j
\right)$, with initial point $q$ and with terminating point $Q$. The
embedding of $\bar{\mathcal{Y}}$ into $\bar{\mathcal{X}} \times
\mathbb{C}$ also defines 
a continuous map of $\gamma+ \Gamma + (-\gamma)$ to
$\mathfrak{P}_n$ which is another geometric braid with base point
$E_n^q$. (The sum of pathes means, that we first walk along $\gamma$
until we reach its terminating point, which is the initial point of
$\Gamma$, then
along $\Gamma$, 
then along $-\gamma$, which means $\gamma$ equipped with orientation
opposite to the orientation of $\gamma$.) After identifying the
pathes $\Gamma_0$ and  $\gamma+ \Gamma + (-\gamma)$ by a
homeomorphism, the two geometric braids are isotopic, since
$\mathcal{X} \setminus \mathcal{X}_0$ does not contain points in the
branch locus of the projection $P_{\mathcal{X}}\mid \mathcal{Y}.\,$
Denote the isotopy class of these braids by $w'_{\mathcal{Y}}$.
Having in mind a homomorphism between $\mathcal{B}_n$ and the group
of isotopy classes of geometric braids with base point $E_n^q$, we
write $w'_{\mathcal{Y}}\overset{\rm def}= w_1^{-1} \cdot
w_{\mathcal{Y}} \cdot w_1$ for a braid $w_1 \in \mathcal{B}_n.\,$

Pick $m$ of the $\mathcal{X}_j,\, j\ge 1,$ and label them in the
order we meet the initial point $q_j$ of $\gamma_j$ when traveling
along $\partial \mathcal{X}_0$ in the direction of orientation of
$\partial \mathcal{X}_0$ as boundary of   $ \mathcal{X}_0$, starting
from $q$. For transparency of the proof we assume $m=k$. The braid
$w$ is quasi-positive. Hence, we may write $w' \overset{\rm def} =
w_1 \cdot w \cdot (w_1)^{-1} = (v_1)^{-1}\, \sigma _{\ell_1}\, v_1
\cdot \ldots \cdot (v_m)^{-1}\, \sigma _{\ell_m}\, v_m $ for $v_j
\in \mathcal{B}_n, \; j=1,\ldots,m.\,$

Consider the $n$ points in $E_n^q$ and label them by
$\xi_1,\ldots,\xi_n$. Put $\varrho = \min\{|\xi_i-\xi_j|: \, i\neq
j\}.\,$ For each $j$ we pick a point $\eta_j \in \mathcal{X}_j$ and
consider the holomorphic curve $\mathcal{Z}_j =\{ (z_1,z_2) \in
\mathcal{X}_j \times \mathbb{C}: \alpha_j\, (z_2
-\xi_{\ell_j})^2=z_1-\eta_j   \}.\,$ Here $\alpha_j$ is a constant
chosen so that $\mathcal{Z}_j $ is contained in $\mathcal{X}_j
\times \{|z_2-\xi_{\ell_j}|<\frac{\varrho}{2} \}.\,$ Consider the
surfaces $\mathcal{Y}_j= \mathcal{Z}_j \cup \bigcup _{i\neq \ell_j,
\ell_{j+1}}\{ (z_1,z_2): z_2=\xi_i,\, z_1 \in \mathcal{X}_j \},\;
j=1,\ldots,m.\,$ Each $\mathcal{Y}_j$ is a Riemann surface which has
$m-1$ connected components and smooth boundary. The boundary
$\partial \mathcal{Y}_j$ corresponds to the conjugacy class of
braids which contains $\sigma _{\ell_j}.\,$ The projection
$P_\mathcal{X} \mid \mathcal{Y}_j:\mathcal{Y}_j \to \mathcal{X}_j$
is a simple branched covering with a single point $\eta_j$ in the
branch locus. The embedding of the closure
$\overline{\mathcal{Y}_j}$ into $\overline{\mathcal{X}_j} \times
\mathbb{C}$ corresponds to a quasi-polynomial, denoted
$\mathcal{P}_j(z,\zeta), \, z \in \bar{\mathcal{X}_j}, \zeta \in
\mathbb{C},\,$ such that the polynomial $\zeta \to
\mathcal{P}_j(\eta_j,\zeta),\, \zeta \in \mathbb{C},\,$ has a double
zero, and the restriction of the quasi-polynomial to
$\bar{\mathcal{X}_j} \setminus \{ \eta_j\}$ is separable.

Let $Q_1$ be the terminating point of $\gamma_1$. Denote by
$\Gamma_1$ a closed arc whose interior parametrizes
$\,\,\partial\mathcal{X}_1 \setminus \{Q_1\}$ with orientation
induced from the orientation of $\partial\mathcal{X}_1$ and whose
endpoints are equal to $Q_1$. Define a mapping from
$\overline{\mathcal{X}_0} \cup \gamma_1 \cup
\overline{\mathcal{X}_1}$ into the space $\overline{\mathfrak{P}_n}$
of all monic polynomials of degree $n$ as follows. For $z_1 \in
\overline{\mathcal{X}_0}$ we let the map be equal to the monic
polynomial with zeros $\overline{\mathcal{Y}} \cap \left(\{ z_1 \}
\times\mathbb{C}\right),\, $ and for $z_1 \in
\overline{\mathcal{X}_1}$ we let the map be equal to the monic
polynomial with zeros $\overline{\mathcal{Y}_1} \cap \left(\{ z_1 \}
\times\mathbb{C}\right).\, $ We extend the map continuously by a map
from $\gamma_1$ to the space $\mathfrak{P}_n$ of polynomials without
multiple roots so that the following holds. The induced mapping from
$\gamma_1 + \Gamma_1 +(-\gamma_1) + \Gamma_0$ to $\mathfrak{P}_n$
with base point $E_n^q$ represents the geometric braid $(v_1)^{-1}\,
\sigma _{\ell_1}\, v_1 \cdot w'_{\mathcal{Y}}$. (We use the same
homomorphism between $\mathcal{B}_n$ and the group of isotopy
classes of geometric braids with base point $E_n^q$ as before.)

Continue by induction. Let $Q_2$ be the terminating point of
$\gamma_2$. 
Denote by $\Gamma_2$ an  arc 
whose interior parametrizes $\,\partial\mathcal{X}_2 \setminus
\{Q_2\}$ with orientation induced from the orientation of
$\partial\mathcal{X}_2$ and whose endpoints are equal to $Q_2$.
Denote by $\Gamma_0^1$ the path obtained by walking along
$\partial\mathcal{X}_0$ according to its orientation from $q=q_1$
until the initial point $q_2$ of the curve $\gamma_2$. Consider the
continuous map from $\overline{\mathcal{X}_0} \cup \gamma_1 \cup
\overline{\mathcal{X}_1} \cup \gamma_2 \cup
\overline{\mathcal{X}_2}$ into the space $\overline{\mathfrak{P}_n}$
which is equal to the previous map on $\overline{\mathcal{X}_0} \cup
\gamma_1 \cup \overline{\mathcal{X}_1}$, which is equal to the monic
polynomial with zeros $\overline{\mathcal{Y}_2} \cap \left( \{ z_2
\} \times\mathbb{C}\right)\, $ for $z_1 \in
\overline{\mathcal{X}_2}$, and which is continuously extended to
$\gamma_2$ so that the following holds. 
The induced mapping from $(\gamma_1 + \Gamma_1 + (-\gamma_1)) +
(\Gamma_0^1 +\gamma_2 + \Gamma_2 +(-\gamma_2) + (-\Gamma_0^1))
+\,\Gamma_0$ has image in $\mathfrak{P}_n$ and represents the braid
$(v_1)^{-1}\, \sigma _{\ell_1}\, v_1 \cdot (v_2)^{-1}\, \sigma
_{\ell_2}\, v_2 \cdot w'_{\mathcal{Y}}.\,$

By induction we obtain a continuous map from
$\mathcal{A}=\overline{\mathcal{X}_0} \cup \bigcup_{j=1,\ldots,m}(
\gamma_j \cup \overline{\mathcal{X}_j})  $ to
$\overline{\mathfrak{P}_n}$ which is holomorphic on
$\bigcup_{j=0,\ldots,m}{\mathcal{X}_j}.\, $ Moreover, the induced
map on the curve $\mathcal{C}$ defines a geometric braid
corresponding to the conjugacy class of $(v_1)^{-1}\, \sigma
_{\ell_1}\, v_1 \cdot \ldots \cdot (v_m)^{-1}\, \sigma _{\ell_m}\,
v_m \cdot w'_{\mathcal{Y}} = w_1 \, w\, w_1^{-1} \,w_1 \,
w_{\mathcal{Y}} \,w_1^{-1} = w_1w w_{\mathcal{Y}}w_1^{-1}   .\,$

The mapping from $\mathcal{A}$ to $\overline{\mathfrak{P}_n}$
corresponds to
a quasi-polynomial of degree $n$ on $\mathcal{A}$ denoted by
$\mathcal {P}$, $\mathcal {P} (z,\zeta) =
\underset{m=0}{\overset{n}{\sum}} \ a_m
(z) \, \zeta^m$, $z \in \mathcal{A}$, $\zeta \in {\mathbb D}.\,$ 
Here $a_n \equiv 1$. The coefficients $a_m$ are continuous functions
on $\mathcal{A}$ which are analytic on the interior ${\rm Int} \,
\mathcal{A} = \underset{j=0}{\overset{k}{\bigcup}} {\mathcal X}_k$.
The polynomial $\zeta \to \mathcal{P}(z,\zeta)$ has multiple zeros
exactly if $z$ is in the branch locus of $P_{\mathcal{X}} \mid
\mathcal{Y}$. This happens for some points in $\mathcal{X}_0$ and
for the points $\eta_j\in \mathcal{X}_j,\,j=1,\ldots,m.\,$ Denote by
$\mathcal{Y}'$ the zero set of the quasi-polynomial in $\mathcal{A}
\times \mathbb{C}.\,$

\smallskip

By the Mergelyan theorem for open Riemann surfaces (\cite{K},
Corollary~4) each coefficient $a_m,\, m=0,1,\ldots,n-1,\,$ can be
uniformly approximated by a holomorphic function $\tilde a_m$ in a
neighbourhood $\tilde {\mathcal{A}}$ of $\mathcal{A}$ on
$\mathcal{X}$. Denote $\tilde{\mathcal P} (z,\zeta) =
\underset{m=0}{\overset{n}{\sum}} \ \tilde a_m (z) \, \zeta^m$, $z
\in \tilde {\mathcal{A}}$, $\zeta \in {\mathbb D}$. Here $\tilde
a_n=a_n \equiv 1.\,$ We may assume that $0$ is a regular value of
$\tilde{\mathcal P}$. If the approximation is good enough then for
each $z \in
\partial {\mathcal{A}}$ the polynomial $\tilde{\mathcal P} (z,\cdot)$ has $n$
distinct roots which are close to the roots of ${\mathcal P} (z,
\cdot)$. There is a thickening $X \subset \tilde {\mathcal{A}}$ of
$\underset{j=0}{\overset{k}{\bigcup}} {\mathcal X}_j \cup
\underset{j=1}{\overset{k}{\bigcup}} \gamma_j$ such that
$\tilde{\mathcal P} (z,\cdot)$ has no multiple roots for $z$ in the
closure of each rectangle $R_j$ added to build the thickening. Put
$\bar Y = \{ (z,\zeta) \in \bar X \times {\mathbb D} :
\tilde{\mathcal P} (z,\zeta) = 0 \}$. The interior $Y$ is
holomorphically embedded into $X \times {\mathbb D}$ and the
projection $P_{\mathcal X} \mid \bar Y: \bar Y \to \bar X$ extends
to a neighbourhood of $\bar Y$ as branched covering with branch
locus in $\underset{j=0}{\overset{k}{\bigcup}} {\mathcal X}_j
\subset X$. In particular, if $\bar X$ is close enough to
$\mathcal{A}$, then the closed geometric braids $\partial Y \subset
\partial X \times {\mathbb C}$ and ${\mathcal Y} \cap ({\mathcal C}
\times {\mathbb C})$ are free isotopic (after identifying $\partial
X$ and $\mathcal{C}$ by a homeomorphism), and, hence, they
correspond to the same conjugacy class in $\mathcal{B}_n$, namely to
the conjugacy class of $w \cdot w_{\mathcal Y} $. The same
conclusion holds for $X$ replaced by any thickening of
$\underset{j=0}{\overset{k}{\bigcup}} {\mathcal X}_j \cup
\underset{j=1}{\overset{k}{\bigcup}} \gamma_j$ contained in $X$. The
lemma is proved. \hfill $\Box$

\medskip

Note that the proof of Lemma 6 provides also a proof of one of the
implications of Rudolph's theorem. Indeed, let $L$ be the closure of
a quasipositive braid. Put $\mathcal{X}=\mathbb{D}$, let the
$\mathcal{X}_j$ be as in the lemma and consider a constant mapping
from $\mathcal{X}_0$ to $\mathfrak{P}_n$. The lemma provides a
simply connected domain $X \subset \mathbb{D}$ and a holomorphic
embedding of a Riemann surface $Y$ into $X \times \mathbb{C}$ such
that $Y\to X$ is a simple branched covering and $Y$ has smooth
boundary $\partial Y$ such that the embedding $\partial Y \to
\partial X \times \mathbb{C}$ defines a closed geometric braid which
corresponds to the conjugacy class of the quasi-positive braid. To
obtain a complex curve in $\mathbb{D} \times \mathbb{C}$ with
boundary contained in $\partial\mathbb{D} \times \mathbb{C}$ and
isotopic to $L$, we map $X$ conformally onto $\mathbb{D}$.

\bigskip

For the proofs of the remaining statements of the theorems we need
the following known proposition.


\begin{prop}
Let $X$ be a connected smooth surface or a connected smooth surface
with boundary. Let $\pi_1 (X,x_0)$ be the fundamental group of $X$
with a given base point $x_0$. The following statements hold.
\begin{enumerate}
\item[1.] There is a one-to-one correspondence between homomorphisms $\Psi : \pi_1 (X,x_0) \to
{\mathcal S}_n$ and  unramified $n$-cove\-rings $p : Y \to X$ with
given label of points in the fiber $p^{-1} (x_0)$.
\item[2.] There is a one-to-one correspondence between homomorphisms $\Phi :
\pi_1 (X,x_0) \to \mathcal{B}_n$ and isotopy classes of separable
quasi-polynomials with fixed value $E_n$ at $x_0$ and given label of
points in $E_n$. The quasi-polynomial lifts a covering $p$ iff the
homomorphism $\Phi$ lifts the homomorphism $p_*:\pi_1 (X,x_0)\to
\mathcal{S}_n $ corresponding to $p$, i.e $p_*=\tau_n \circ \Phi$
for the canonical homomorphism $\tau_n:\mathcal{B}_n \to
\mathcal{S}_n $.
\item[3.] The connected components of the covering space $Y$ of an unramified holomorphic
$n$-covering $p:Y\to X$ are in bijective correspondence to the
orbits of $p_* (\pi_1 (x,z_0)) \subset {\mathcal S}_n$ on the set
consisting of $n$ points. In particular, $Y$ is connected iff $p_*
(\pi_1 (X,x_0))$
acts transitively. 
\item[4.] Suppose $X$ is a $2$-manifold of genus $g$ with boundary. Suppose the boundary $\partial X$ is connected and the base point $x_0$ is
contained in the boundary. Denote by $\{\partial X\}$ the element of
the fundamental group $\pi_1(X,x_0)$ which is represented by
traveling along $\partial X$ in the direction of orientation as
boundary of $X$ starting from $x_0$. (So, $\{\partial X\}$ is the
product of $g$ commutators of suitable generators of the fundamental
group.) Let $p:Y \to X$ be an unramified $n$-covering and let
$\mathcal{P}(x,\zeta)$ be a separable quasi-polynomial lifting it,
i.e. the covering is equivalent to $P_X \mid
\mathfrak{S}_{\mathcal{P}}: \mathfrak{S}_ {\mathcal{P}}\to X$.

Then the connected components of the boundary $\partial \, Y$
correspond to the orbits of the single even permutation $p_* (\{
\partial X \})$. $p_* (\{ \partial X \})$ is the product of $g$
commutators in ${\mathcal S}_n$. 

Moreover, the free isotopy class of the boundary link $\,\partial
\mathfrak{S}_{\mathcal{P}} \subset \partial X \times \mathbb{C}$ is
a closed geometric braid representing the conjugacy class of the
product of $g$ commutators in $\mathcal{B}_n$.
\end{enumerate}
\end{prop}

Notice that a conjugate of a commutator is again a commutator.
\bigskip

\section{Proof of the remaining statements of the theorems}

\noindent \textbf{Proof of statement 3 of Theorem 1.}\\
\noindent \textbf{The sharpness of the estimate for the general case
of links $L$ (statement 3)} is an immediate consequence of the
following lemma.

\begin{lm}Let $\mathcal {X}$ be an open Riemann surface with smooth
connected boundary of genus $g(\mathcal{X})\ge 1$. For any natural
$n \ge 1$ there exists an unbranched holomorphic $n$-covering
$\mathcal {Y}$ of $\mathcal {X}$ such that $\mathcal {Y}$ is
connected and has $n$ boundary components. Moreover, there is a
holomorphic embedding of $\mathcal {Y}$ into the disc bundle
$\mathcal {X} \times \mathbb {D}$ that lifts the
covering map.
\end{lm}

Indeed, let $K$ be a smoothly analytic knot, i.e. $K=\partial
\mathbb{B}^2 \cap \mathcal{X}$ for a complex curve ${\mathcal X} =
\{ z$ in a neighbourhood of $\overline{{\mathbb B}^2} : f(z) = 0
\}$. Here $f$ is an analytic function with non-vanishing gradient in
a neighbourhood of $\mathcal{X}$. Suppose $g(\mathcal{X}) \ge 1$. We
may assume that $\mathcal{X}\cap \mathbb{B}^2$ is diffeomorphic to
$\mathcal{X}$. Identify the disc bundle $\mathcal{X} \times
\mathbb{D}$ with a small tubular neighbourhood of ${\mathcal X}$ in
${\mathbb C}^2$ and consider the Riemann surface $\mathcal{Y}$ of
Lemma 7 to be an embbedded submanifold of the tubular neighbourhood.
For these $\mathcal{X}$ and $\mathcal{Y}$ equality in the
Riemann-Hurwitz relation \eqref{eq4'} 
is obtained, since the covering is unbranched, ${\mathcal Y}$ is
connected and the number of boundary components of ${\mathcal Y}$ is
maximal. Let $L=\mathcal{Y}\cap \partial \mathbb{B}^2$.
If the tubular neighourhood is small enough then $\mathcal{Y}$ is
diffeomorphic to $\mathcal{Y} \cap \mathbb{B}^2$. It follows that in
this case equality is attained in \eqref{eq3}. 
Statement 3 is proved.

\medskip

\noindent {\bf Proof of Lemma 7.} Let $\mathcal{F}:\mathbb{C} \to
\mathbb{C}$ be the complex linear mapping $\mathcal{F}(\zeta)=
e^{\frac{2\pi i}{n}} \,\zeta.\,$ The $n$-th iterate $\mathcal{F}^n$
is the identity. The fundamental group $\pi_1(\mathcal{X},x_0)$ is a
free group on $g=g(\mathcal{X})$ generators. Denote by $\langle
\mathcal{F} \rangle$ the group of self-homeomorphisms of
$\mathbb{C}$ generated by $\mathcal{F}$. Let
$\Phi:\pi_1(\mathcal{X},x_0) \to \langle \mathcal{F} \rangle $ be a
homomorphism which assigns the element $\mathcal{F}$ to one of the
generators of the fundamental group, and assigns to each of the
other generators any element of the group $\langle \mathcal{F}
\rangle$. By proposition 4 the image $\Phi(\{\partial
\mathcal{X}\})$ is the identity.

Let $\tilde{\mathcal{X}}$ be the universal covering of $X$. The
fundamental group $\pi_1(\mathcal{X},x_0)$ acts on
$\tilde{\mathcal{X}}  \times \mathbb{C}$ as follows.
\begin{equation}\label{eq5}
\tilde {\mathcal{X}} \times \mathbb{C} \ni (x,\zeta) \to (\gamma(x),
\Phi(\gamma)(\zeta)),\;\; \gamma \in \pi_1(\mathcal{X},x_0) \,.
\end{equation}
The action is free and properly discontinuous. Hence, the mapping
\begin{equation}\label{eq6}
p:\tilde{\mathcal{X}} \times \mathbb{C} \to  \mathcal{E}\overset{\rm
def}{=}\tilde{\mathcal{X}} \times \mathbb{C} \diagup
\pi_1(\mathcal{X},x_0)
\end{equation}
is a holomorphic covering map. It defines a holomorphic fiber bundle
over $\mathcal{X}$ with fiber being a complex line, and with
transition functions being complex linear. Since $\mathcal{X}$ is
open this holomorphic line bundle is trivial. The covering map
\eqref{eq6} respects a holomorphic foliation on $\tilde{\mathcal{X}}
\times \mathbb{C}$, namely the trivial foliation with leaves
$\tilde{\mathcal{X}}\times \{\zeta\},\, \zeta \in \mathbb{C}.\,$ The
map $\mathcal{F}$ permutes the points of the set $E=\{
1,e^{\frac{2\pi i}{n}},\ldots, e^{\frac{2\pi i (n-1)}{n}}\} $ along
a cycle of length $n$. Hence, the covering map \eqref{eq6} maps the
set $\tilde{\mathcal{X}} \times E$ to a single leaf. This leaf is a
Riemann surface which we denote by $\mathcal{Y}$ and identify with a
surface which is $n$-horizontally embedded into the trivial bundle
$\mathcal{X}\times \mathbb{C}$. The projection
$P_{\mathcal{X}}:\mathcal{X} \times \mathbb{C}  \to \mathcal{X} $
restricts to $\mathcal{Y}$ as unramified covering. The covering
corresponds to a homomorphism $\Psi: \pi_1(\mathcal{X},x_0) \to
\langle s \rangle  $, where $s=(12\ldots n)$ is a cycle of length
$n$, and $\langle s \rangle$ is the subgroup of the symmetric group
generated by $s$. Hence, $\Psi( \{\partial \mathcal{X}\})= \rm{id}$,
and $Y$ has $n$ boundary components by Proposition 5. \hfill $\Box$

\medskip

\noindent \textbf{Proof of statement 4 of Theorem 1.}

\noindent \textbf{The sharpness of the bound for the case when $L$
is also required to be a knot (statement 4)} is obtained as follows.
Let $K$ be an analytic knot with $g_4 (K) = g \ge 1$. Using the
isotopy
provided by Lemma 4 we may assume that 
$K = {\mathcal X} \cap
\partial {\mathbb B}^2$ for a smooth complex curve ${\mathcal X}$
contained in 
$(1+\epsilon)\mathbb{D} \times \epsilon \mathbb{D}$ for a small
positive number $\epsilon$  such that $P_1 \mid {\mathcal X} :
{\mathcal X} \to (1+\varepsilon) {\mathbb D}$ is a branched covering
with branch locus in ${\mathbb D}$ and $\mathcal{X}$ is
diffeomorphic to $\mathcal{X} \cap \mathbb{B}^2$.

\smallskip

{\bf If $n$ is odd} there exists an unbranched holomorphic covering
$p : {\mathcal Y} \to {\mathcal X}$ such that ${\mathcal Y}$ has
connected boundary. This follows immediately from Proposition 5 and
a theorem of Ore which says that each even permutation is a
commutator. For convenience of the reader we provide the following
simple examples on commutators. Example 1 together with statements 1
and 4 of Proposition 5 provide the required unbranched covering.
Indeed, take the covering corresponding to the homomorphism $\Psi:
\pi_1(\mathcal{X},x_0) \to \mathcal{S}_n$ for which
$\Psi(a_1)=s_1,\, \Psi(a_2)=s_2,\, \psi(a_j)={\rm id}, \,
j=3,\ldots,2g,\,$ for a suitable choice of generators $a_j$ of
$\pi_1(\mathcal{X},x_0)$ and for the
 permutations $s_1$ and $s_2$ of Example 1. \\
Example 2 will be used below.

\bigskip

\noindent {\bf Example 1.} Suppose $n$ is an odd number, $n=2m+1$.
Consider the following two permutations $s_1 = (23) \ldots (2m \
2m+1)$ and $s_2 = (12)(34) \ldots (2m-1 \ 2m)$ in $\mathcal{S}_n$.
Then the commutator $s = [s_1 , s_2]$ is a cycle of order $n$. (See
fig.~6a for $n=7$.)

\medskip

\noindent {\bf Example 2.} Let $n=2m$ be an even number. Consider
the permutations $s_1 = (23) \ldots (2m - 2 \ 2m-1)$ and $s_2 = (12)
\ldots (2m-1 \ 2m)$ in $\mathcal{S}_n$. The commutator $s = [s_1 ,
s_2]$ is the disjoint union of two cycles of order $\frac{n}{2}$.
Note that the subgroup $\langle s_1 , s_2 \rangle$ of
$\mathcal{S}_n$ generated by $s_1$ and $s_2$ acts transitively on
$\{1,2,\ldots,n\}$. (See fig.~6b for $n=8$.)

$$
\includegraphics[width=10cm]{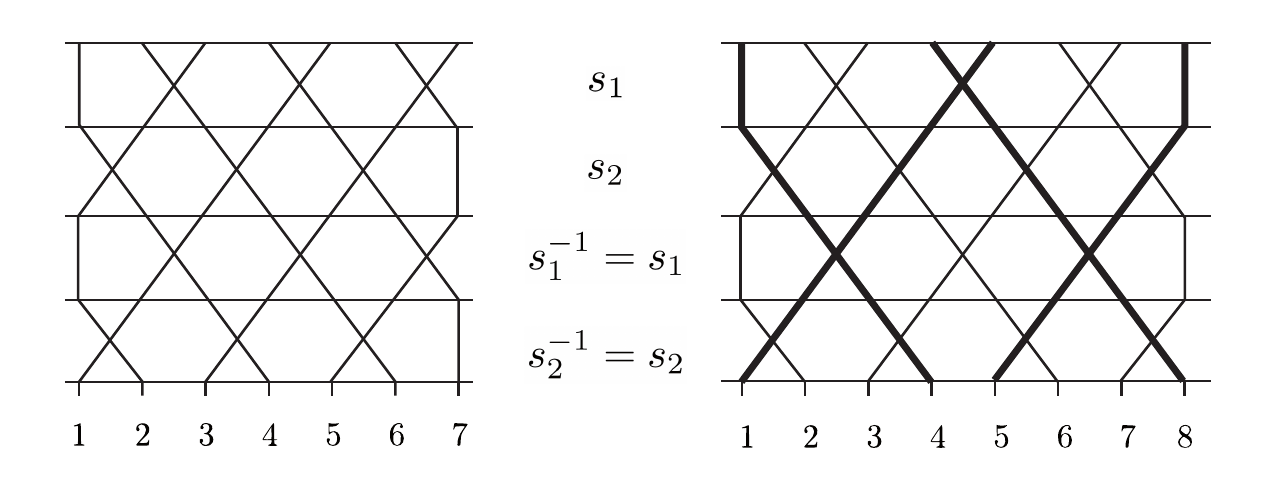}
$$
$$\mbox{fig. 6a} \qquad\qquad\qquad\qquad\qquad\qquad\qquad\qquad \mbox{fig. 6b}$$

\centerline{Figure 6.}

\bigskip

Take the complex structure on $\mathcal{Y}$ which is induced by the
covering map $p$. Consider the holomorphic mapping $i = (p,f) :
{\mathcal Y} \to {\mathcal X} \times  \, {\mathbb D}$. Here $f$ is a
bounded holomorphic function. 
The mapping $i$ is an immersion. It is an embedding of $p^{-1}
({\mathcal X}_0)$ into ${\mathcal X}_0 \times \,{\mathbb D}$ for a
domain ${\mathcal X}_0 \subset {\mathcal X}$ iff $f$ separates the
points of the fiber $p^{-1} (z)$ for each $z \in {\mathcal X}_0$.
Choose the holomorphic function $f:\mathcal{Y} \to  \mathbb{D}\,$ so
that it separates all points of the fibers $p^{-1} (z_j)$ for all
points $z_j$ in the branch locus of $P_1 \mid {\mathcal X}$. This is
possible since ${\mathcal X}$ is an open Riemann surface with smooth
connected boundary, hence it has bounded holomorphic functions with
prescribed values at finitely many points. From the identity theorem
applied to $f \circ \varphi_j$ for the local inverses $\varphi_j$ of
$p$ it follows that the function $f$ separates points of fibers for
all $z \in {\mathcal X}$ not belonging to a discrete subset
$\Lambda$ of ${\mathcal X}$. The disc $\bar{\mathbb D}$ contains
only finitely many points of $P_1 (\Lambda)$. Join them by pairwise
disjoint simple smooth arcs with points in $\partial {\mathbb D}$.
Remove small neighbourhoods of the arcs from $\bar{\mathbb D}$ so
that we obtain a simply connected smoothly bounded domain ${\mathcal
D} \subset {\mathbb D}$. We may assume that ${\mathcal D}$ contains
the branch locus of $P_1 \mid {\mathcal X}$, in other words,
${\mathcal X}_1 = {\mathcal X} \cap ({\mathcal D} \times {\mathbb
C})$ is diffeomorphic to
${\mathcal X}.\,\,$
Identify the disc bundle ${\mathcal X}_1 \times  \, {\mathbb D}$
with a small enough tubular neighbourhood of ${\mathcal X}_1$ in
$\mathbb{C}^2$. We obtain an embedding of ${\mathcal Y}_1 = p^{-1}
({\mathcal X}_1)$ into the tubular neighbourhood of ${\mathcal X}_1$
such that ${\mathcal Y}_1$ is diffeomorphic to ${\mathcal Y}$. By
the Riemann-Hurwitz relation we obtain
$$
1-2g({\mathcal Y_1}) = \chi ({\mathcal Y}_1) = n \cdot \chi
({\mathcal X}_1) = n (1-2g({\mathcal X}_1)) \, ,
$$
hence, since $\mathcal{X}_1$ is diffeomorphic to  $\mathcal{X}$ and
${\mathcal Y}_1$ is diffeomorphic to ${\mathcal Y}$, we obtain
$$
g({\mathcal Y}) = n \cdot g ({\mathcal X}) - \frac{n-1}{2} \, .
$$
Lemma~5 provides a further isotopy to an analytic knot $K$ and an
analytic knot $L$ which is an $n$-braided link around $K$ contained
in a small tubular neighbourhood of $K$ such that equality in
relation \eqref{eq4} for the $4$-ball genus bound is attained.

\medskip

It remains to consider the case when {\bf $n$ is even}. We may
assume that $\mathcal{X}$ is as in the case of odd $n$. There exists
an unbranched holomorphic $n$-covering $p : {\mathcal Y}\to{\mathcal
X}$ with ${\mathcal Y}$ connected and with boundary consisting of
two connected components. This follows from statements 1 and 4 of
Proposition 5 and Example 2 above. (See fig. 6b for the case $n=8$.)
By statement 2 of Proposition 5 there is a smooth separable
quasipolynomial $\mathcal{P}$ that lifts the covering. Hence, we may
identify $\mathcal{Y}=\mathfrak{S}_{\mathcal{P}}$ with the zero set
of the quasipolynomial so that the covering map $p$ equals
$P_{\mathcal{X}} \mid \mathcal{Y}: \mathcal{Y \to \mathcal{X}}.\,$
Shrinking $\mathcal{X}$ we may assume that $\mathcal{X}$ and
$\mathcal{Y}$ have smooth boundary and $\bar{\mathcal{Y}}$ is
smoothly $n$-horizontally embedded into $\bar{\mathcal{X}} \times
\mathbb{D}.\,$


Let $w_{{\mathcal Y}} \in {\mathcal B}_n$ represent the conjugacy
class corresponding to the isotopy class of $\partial {\mathcal Y}$
in $\partial {\mathcal{X}} \times \, \mathbb D$. Take for $w$ a
conjugate of a generator of ${\mathcal B}_n$ which permutes two
strands of $w_{{\mathcal Y}}$ corresponding to different connected
components of $\partial {\mathcal Y}$ (in other words, the closure
of $w \cdot w_{{\mathcal Y}}$ defines a connected closed braid).
Apply statement 2 of Proposition 4. We obtain a domain $X_1 \subset
{\mathcal X}$ of the form $X_1 = P_1^{-1} ({\mathcal D}_1)$ for a
smoothly bounded simply connected domain ${\mathcal D}_1,\,$
${\mathcal D}_1 \subset \mathbb{D}$, such that $X_1$ is
diffeomorphic to ${\mathcal X}$. Also we obtain an embedding of an
open Riemann surface $i : Y_1 \to X_1 \times  \, {\mathbb D}$ into
the disc bundle such that $P_{\mathcal X} \mid Y_1$ is a branched
holomorphic covering with a single branch point and with connected
boundary $\partial Y_1 \subset
\partial X_1 \times  \, {\mathbb D}$ which determines a closed geometric braid
which (after identifying $\partial \mathcal{X} $ with $\partial X_1$
by a homeomorphism) represents the conjugacy class of $w \cdot
w_{{\mathcal Y}}$.

By the Riemann-Hurwitz relation
$$
1-2g (Y_1) = n(1-2g (X_1)) - 1 \, ,
$$
hence
$$
g (Y_1) = ng (X_1) - \frac{n-2}{2} \, .
$$

Since $X_1$ is diffeomorphic to $\mathcal{X}$, and, hence, to
$\mathcal{X}\cap \partial \mathbb{B}^2$, we have $g_4(K)=g(X_1)$.
Lemma~5 gives a further isotopy of $K$ to a smoothly analytic knot
again denoted by $K$ and an analytic knot $L$ in a small
neighbourhood of $K$ which is an $n$-braided link around $K$ such
that $g_4 (L) = g(Y_1)$. The first part of statement 4 is proved.

\medskip

The following example proves the last part of statement 4. Embed the
standard punctured torus holomorphically into $\mathbb{C}^2$ using
the Weierstra{\ss} $\wp$-function:
\begin{equation}\label{eq7}
\left(\mathbb{C} \setminus (\mathbb{Z}+i\mathbb{Z})\right)\diagup
(\mathbb{Z}+i\mathbb{Z}) \ni \zeta \to (\wp(\zeta), \wp'(\zeta)) \in
\mathbb{C}^2\,.
\end{equation}
Denote the image of the embedding by $\mathcal{X}$. Let $R$ be a
large positive number. The intersection $X_R= \frac{1}{R}\mathcal{X}
\cap  \mathbb{B}^2$ is a torus with a hole. If $R$ is large then
$X_R$ contains a domain $\mathcal{R}$ which is adjacent to $\partial
X_R$ and is conformally equivalent to an annulus of conformal module
larger than $\frac{\pi}{2} \,\frac{1}{\log \frac{3+\sqrt{5}}{2}}$.
(Recall that for $0\leq r_1<r_2 \leq \infty$ the conformal module of
the annulus $\{r_1<|z|<r_2\}$ in the complex plane equals
$\frac{1}{2\pi} \log\frac{r_2}{r_1}.\,$) Put $X= X_R \subset
\mathbb{B}^2$ for a number $R$ with this property. Let $K=\partial X
\subset
\partial \mathbb{B}^2.\,$ Then $K$ is a smoothly analytic knot.

Suppose for any $a>0$ there exists an analytic knot $L$ contained in
the tubular neighbourhood $N(K)= \partial \mathbb{B}^2 \cap
\mathcal{T}_a(\frac{1}{R}\mathcal{X})$, such that $n=w_{N(K)}(L)=3$
and equality is obtained in the $4$-ball genus estimate \eqref{eq4}
for $K$ and $L$. Let $Y$ be the complex curve bounded by $L$. Apply
Lemma 2 and Proposition 1. 
Let $\mathcal{H}$ be the Levi-flat hypersurface of Proposition 1 and
let $A$ be the set defined before the statement of Proposition 1.
Put $X'=X \setminus \overline A\,,$ $Y'=Y \setminus \overline A\,,$
$L'=Y \cap \mathcal{H}\,$. If $a>0$ is small then $A$ is contained
in a small neighbourhood of $K$ in $\mathbb{C}^2$.

If $R$ is large then the set $\mathcal{R}' \overset {\rm def} = X'
\cap \mathcal{R}$ is conformally equivalent to an annulus of
conformal module close to that of $\mathcal{R}$, in particular the
conformal module of $\mathcal{R}'$ is larger than $\frac{\pi}{2}
\,\frac{1}{\log
\frac{3+\sqrt{5}}{2}}$. 

Apply Proposition 2 with $\mathcal{X}=X'$ and $\mathcal{Y}=Y'$:
Since $n=3$ is odd, $L=\partial Y$ is connected and equality holds
in \eqref{eq4}, by Proposition 2 the Riemann surface
$\mathcal{Y}=Y'$ has connected boundary and the covering is
unramified. 
The embedding of $Y'$ in the disc bundle over $X'$
defines a holomorphic map of $X'$ to $\mathfrak{P}_n$. Its
restriction to $\mathcal{R}'$ is a holomorphic map
into $\mathfrak{P}_n$ which represents the free isotopy class of
$L'=\partial Y'$ (a commutator class by Proposition 4). Since the
conformal module of $\mathcal{R}'$ is large, by Lemma 8.3 of
\cite{J2} the class of $L'$ is the conjugacy class of a pure braid,
i.e. $L'$ cannot be connected. The contradiction proves the last
part of statement 4 of Theorem 1. Theorem~1 is proved. \hfill $\Box$

\bigskip

\noindent {\bf Proof of Theorem 2.} Statement~1 follows from
Proposition~$1'$.

We will prove now statement 2. Let $\tilde X$ and $\tilde Y$ be
relatively closed complex curves in a neighbourhood $\tilde \Omega$ of $\bar{\mathbb{B}}^2$ such that $K=\tilde X \cap
\partial \mathbb{B}^2$ and $L=  \tilde Y \cap \partial \mathbb{B}^2$. For the domain $\Omega_1$ of Proposition $1'$ (with $\Omega =\mathbb{B}^2$) we let $X_1= \tilde X \cap \Omega
_1$, $Y_1 = \tilde X \cap \Omega_1$, $K_1=\partial X_1$ and $L_1 =
\partial Y_1$.
Consider the bundle $\mathcal{T}(\overline{X_1})\to \overline{X_1}$
of proposition 3 with trivializaton inducing Seifert framing for
$K_1$ on $\partial \Omega_1$. The intersection of $Y_1$ with the
discs of the trivialized bundle defines a continuous mapping from
$X_1$ to the space $\overline{\mathfrak{P}_n}$ of monic polynomials
of degree $n$. The restriction of this mapping to $\partial X_1$ is
a mapping to the space $\mathfrak{P}_n$ of monic polynomials of
degree $n$ without multiple zeros which represents the pattern of
the braided link $\partial Y_1$ .

If the covering $p:\overline{Y_1}\to \overline{X_1}$ induced by the
bundle projection is unramified then by Proposition 5 the pattern
$\mathcal{L}_1$ of the link $\partial Y_1$ is the conjugacy class of
a product of $g=g(X_1)=g_4(K)$ commutators in $\mathcal{B}_n$.




Consider now the general case. Let $d \subset X_1$ be a smoothly
bounded simply connected domain which contains the branch locus of
$p$ such that $\bar{X_1} \backslash d$ is diffeomorphic to
$\bar{X_1}$ (in particular, $\partial {X_1} \cap
\partial d \ne \emptyset$). (For example, one can take for $d$ the union of
the following sets: suitable neighbourhoods of simple disjoint arcs
joining a critical value of $p$ with a boundary point of $X_1$, and
a suitable simply connected part of a collar of $\partial X_1$ in
$X_1$.) Put $\mathcal{Y}_d \overset{\rm def}{=} p^{-1} (d)$ and use
the following notation: $\mathcal{X}_{d}\overset{\rm def}{=} d$,
$\mathcal{X}_{Cd}\overset{\rm def}{=}{X_1} \backslash \bar d$ and
${\mathcal Y}_{Cd} \overset{\rm def}{=} p^{-1} (\mathcal{X}_{Cd})$.

Take a base point $q$ in $\bar{X_1}$ which is a boundary point of
$X_1 \cap \partial d$ (in particular, $q \in
\partial d \cap \partial {X_1}$). Choose the point $E_n^q =
P_{\mathbb D}(\mathcal{Y} \cap (\{q\} \times \mathbb D))$ as base
point in the symmetrized configuration space.

Let $\Gamma_d$ be an arc  whose interior parametrizes
$\partial \mathcal{X}_{d} \setminus \{q\}$ 
and whose two endpoints are equal to $q$. Respectively, we denote by
$\Gamma_{Cd}$ an arc whose interior
parametrizes $\partial \mathcal{X}_{Cd} \setminus \{q\}$ 
and whose two endpoints are equal to $q$. Both arcs are equipped
with orientation induced by orienting $\partial \mathcal{X}_{d}$, or
$\partial\mathcal{X}_{Cd}$, respectively, as boundaries of the
domains $\mathcal{X}_{d}$, and $\mathcal{X}_{Cd}$, respectively.

The $n$-horizontal embeddings $\partial\mathcal{Y}_{d} \subset
\partial\mathcal{X}_{d} \times \mathbb{C} $, and  $\partial\mathcal{Y}_{Cd} \subset
\partial\mathcal{X}_{Cd} \times \mathbb{C} $ respectively, define
continuous mappings from $\Gamma_d$ into $\mathfrak{P}_n$, and from
$\Gamma_{Cd}$ into $\mathfrak{P}_n$, respectively. Identifying
$\mathcal{B}_n$ with the group of isotopy classes of geometric
braids with base point $E_n^q$ we obtain (after identifying
$\Gamma_{d}$ and $\Gamma_{Cd}$ with the unit interval) two braids
$w_d$ and $w_{Cd}$. Since $d$ is simply connected the braid $w_d$ is
quasi-positive by Rudolph's theorem. Since the covering over
$\mathcal{X}_{Cd}$ is unramified the braid $w_{Cd}$ is a product of
$g$ commutators in $\mathcal{B}_n$. (See statement 4 of Proposition
5.)

Let $\Gamma_{X_1}$ be an arc 
whose interior parametrizes $\partial X_1 \setminus \{q\}$ (with
orientation induced from the orientation of $\partial X_1$ as
boundary of $X_1$) and whose two endpoints are equal to $q$. The
isotopy class of the continuous mapping from $\Gamma_{X_1}$ to
$\mathfrak{P}_n$ which is defined by the $n$-horizontal embedding
$\partial Y_1 \subset
\partial X_1 \times \mathbb{C} $ is equal (after identification
of the curves $\Gamma_{d}+\Gamma_{Cd}$ and $\Gamma_{X_1}$) to the
braid $w_d \cdot w_{Cd}$. Statement 2 is proved.

\smallskip

It remains to prove statement~3. By Lemma 4 after an isotopy we are
in the situation when the knot is equal to 
${\mathcal X} \cap
\partial {\mathbb B}^2$ for a smooth relatively closed complex curve ${\mathcal X}$ in
a small neighbourhood of $\{ z_2 = 0 \}$ such that for a small
positive number $\epsilon$ the mapping $P_1 : {\mathcal X} \to
(1+\varepsilon) \, {\mathbb D}$ is a branched covering with branch
locus in ${\mathbb D}$ and $\mathcal{X}\cap \mathbb{B}^2$ is
diffeomorphic to $\mathcal{X}$. We may assume that $\mathcal{X}$ has
smooth boundary.

\smallskip

Let $[\alpha_1 , \beta_1] \cdot \ldots \cdot [\alpha_g , \beta_g]$,
$\alpha_j , \beta_j \in {\mathcal B}_n$ for $j = 1,\ldots , g$, be a
product of $g=g(\mathcal{X})$ commutators. By statements 2 and 4 of
Proposition~5 there exists 
a smooth $n$-horizontal embedding of a smooth surface with boundary
$\overline{\mathcal{Y}}$, $\overline{\mathcal{Y}} \to
\overline{\mathcal{X}} \times \mathbb{C},\,$ such that the embedding
of the boundary $\partial{\mathcal{Y}} \to \partial {\mathcal{X}}
\times \mathbb{C}\,$ corresponds to the conjugacy class of the afore
mentioned product of commutators. Indeed, let
$\overline{\mathcal{Y}}$ be the zero set of the quasi-polynomial
which corresponds to the homomorphism $\Phi$ for which
$\Phi(a_j)=\alpha_j,\, \Phi(b_j)=\beta_j\,$ for suitable generators
$a_j$ and $b_j,\, j=1,\ldots,g,\,$ of the fundamental group.

Let $w \in \mathcal{B}_n$ be a quasipositive braid such that the
pattern $\mathcal{L}$ is the conjugacy class of the braid $w \cdot
[\alpha_1 , \beta_1] \cdot \ldots \cdot [\alpha_g , \beta_g]$. By
Proposition 4 there is an open subset $X_1$ of $\mathcal{X}$
diffeomorphic to $\mathcal{X}$ and a holomorphically embedded
manifold $Y_1 \subset X_1 \times \mathbb{C}$ so that for the
canonical projection $P_{\mathcal{X}}$ the mapping
$P_{\mathcal{X}}|Y_1:Y_1 \to X_1$ is a holomorphic branched
$n$-covering. The number of branch points $B$ equals the exponent
sum of the braid $w$. Moreover, the isotopy class of the link
$\partial Y_1 \subset
\partial X_1 \times \mathbb{C}$ corresponds to the conjugacy class
of the braid $w \cdot [\alpha_1 , \beta_1] \cdot \ldots \cdot
[\alpha_g , \beta_g]\,.$

By Lemma 5 the conjugacy class $\mathcal{L}$ of $w \cdot [\alpha_1 ,
\beta_1] \cdot \ldots \cdot [\alpha_g , \beta_g]\,,\,$ can be
realized by an analytic link contained in an a priory given
neighbourhood of a knot $K \subset \partial \mathbb{B}^2$ which is
isotopic to $\mathcal{X}\cap \partial \mathbb{B}^2$.


\smallskip

Theorem 2 is proved. \hfill $\Box$

\bigskip
\noindent {\bf Proof of Lemma 1.} The proof uses the proof of
statement 3 of Theorem 2. We may assume after an isotopy which moves
$K$ and $L$ that the knot $K$ has the form $\mathcal{X} \cap
\partial \mathbb{B}^2$ with $\mathcal{X}$ as in the proof of
statement 3. The pattern $\mathcal{L}$ of $L$ is the closure of a
quasipositive braid $w$ (it corresponds to the case of Theorem 2
when the product of commutators is the identy). Let $B$ be the
exponent sum of $w$. For $Y_1$ and $X_1$ as in the proof of of
statement 3 of Theorem 2 the following equation for the Euler
characteristic holds
\begin{equation}\label{eq8}
\chi(Y_1) = n \chi(X_1)-B \,.
\end{equation}
Under the condition of Lemma 1 the closure of the quasipositive
braid $w$ is connected. Notice the following fact. For a knot $L'$
in a tubular neighbourhood of the unknot
which represents the closure of a quasipositive braid with exponent
sum $B$ we have
\begin{equation}\label{eq9}
1-2g_4 (L') =  n-B\,.
\end{equation}
Indeed, we may assume after an isotopy that  $L' =
\partial {\mathbb B}^2 \cap Y'$ for a complex curve $Y'$ in a
neighbourhood of $\bar{\mathbb B}^2$ such that $P_1 \mid Y' : Y' \to
(1+\epsilon){\mathbb D}$ is a branched $n$-covering with branch
locus in ${\mathbb D}$ and $Y'$ is diffeomorhic to $Y' \cap
\mathbb{B}^2$. Then $g_4 (L') = g(Y')$. Hence,
$$
1-2g_4 (L') = \chi (Y') = n
\chi\left((1+\epsilon)\mathbb{D}\right)-B= n-B\,.
$$
Equations \eqref{eq8} and \eqref{eq9} with the pattern $\mathcal{L}$
of Lemma 1 instead of $L'$ imply
$$
1-2g(Y_1)=n(1-2g(X_1))-B=-2ng(X_1) +(n-B)=-2ng(X_1)
+1-2g_4(\mathcal{L}).
$$
As in the proof of statement 3 of Theorem 2 Lemma~5 provides an
isotopy of $\mathcal{X} \cap \partial \mathbb{B}^2$ to a smoothly
analytic knot $\tilde K$, which bounds a Riemann surface $\tilde X
\subset {\mathbb B}^2$ that is diffeomorphic to ${\mathcal X}$ and
thus to $X_1$, and an $n$-braided link $\tilde L$ in a small tubular
neighbourhood of $\tilde K$, which bounds a Riemann surface $\tilde
Y$ that is diffeomorphic to $Y_1$. Moreover, the pattern of $\tilde
L$ is $\mathcal{L}$. Then $g_4(\tilde K)=g_4(K)$ and $g_4(\tilde
L)=g_4(L)$ (since the isotopy which moves $K$ to $\tilde K$ moves
$L$ to a knot  in a tubular neighbourhood of $\tilde K$ which is
isotopic to $L$ and has the same pattern as $L$ and thus as $\tilde
L$).
Hence,
$$
g_4(L)=g(\tilde Y)=g(Y_1)=ng(X_1) +g_4(L)= ng_4(K)+g_4(\mathcal{L})
$$
The lemma is proved. \hfill $\Box$

\bigskip

\noindent \textbf{Added in  proof.} The proof of Theorem 1 gives immediately also the following variant of the first two statements for the $4$-ball Euler characteristics $\chi _4$ of links (defined similarly as the $4$-ball genus):

{\it For each smoothly analytic knot $K \subset \partial \mathbb{B}^2$ there is a tubular neighbourhood $N(K)\subset \partial \mathbb{B}^2$ of $K$ such that for each analytic link $L \subset N(K)$ with winding number $n \ge 1$ the inequality $\chi_4(L) \leq n \, \chi_4(K)$ holds. }

Indeed, notice that statement (15) is also true if $Y$ is not connected. One has to choose the algebraic curve $C$ in the proof of statements (1) and (2) to be connected.

\bigskip

\setcounter{lm}{0} \setcounter{prop}{0} \setcounter{cor}{0}

\end{document}